\documentclass{amsart}
\usepackage{amssymb,amsmath,amsfonts,epsfig,latexsym,amsthm}

\newtheorem{thm}{Theorem}[section]
\newtheorem{defn}[thm]{Definition}
\newtheorem{prop}[thm]{Proposition}
\newtheorem{lem}[thm]{Lemma}

\newtheorem*{TQ}{Totaro's Question}
\newtheorem*{KCT}{Klyachko's Classification Theorem}

\theoremstyle{definition}
\newtheorem{remk}[thm]{Remark}
\newtheorem{prob}[thm]{Problem}
\newtheorem{example}[thm]{Example}

\def\A{\ensuremath{\mathbf{A}}}
\def\C{\ensuremath{\mathbf{C}}}
\def\P{\ensuremath{\mathbf{P}}}
\def\Q{\ensuremath{\mathbf{Q}}} 
\def\R{\ensuremath{\mathbf{R}}}
\def\Z{\ensuremath{\mathbf{Z}}}

\def\E{\mathcal{E}}
\def\F{\mathcal{F}}
\def\L{\mathcal{L}}

\def\O{\ensuremath{\mathcal{O}}}

\def\Po{\mathcal{P}}

\def\DDelta{\mathbf{\Delta}}
\def\Puck{\mathbf{Puck}}

\def\SSigma{\mathbf{\Sigma}}
\def\Star{\mathbf{Star}}

\def\CX{\mathbf{CX}}
\def\CY{\mathbf{CY}}
\def\bf{\mathbf{f}}
\def\bg{\mathbf{g}}
\def\bu{\mathbf{u}}
\def\bw{\mathbf{w}}

\def\FFl{\mathcal{F}\ell}
\def\Fl{F\ell}
\def\PL{P\!L}
\def\PP{P\!P}

\DeclareMathOperator{\divisor}{div}
\DeclareMathOperator{\Hom}{Hom}
\DeclareMathOperator{\relint}{rel \, int}

\DeclareMathOperator{\Sym}{Sym}
\DeclareMathOperator{\Tr}{Tr}

\def\<{\ensuremath{\langle}}
\def\>{\ensuremath{\rangle}}

\begin{document}
    
\title[Branched covers of fans and the resolution property]{Toric vector bundles, branched covers of fans, and the resolution property}           

\subjclass[2000]{}
    
\author{Sam Payne}
\address{Stanford University, Mathematics, Bldg. 380, Stanford, CA 94305}
\email{spayne@stanford.edu}
\thanks{Supported by a Graduate Research Fellowship from the NSF}

\begin{abstract}
We associate to each toric vector bundle on a toric variety $X(\Delta)$ a ``branched cover" of the fan $\Delta$ together with a piecewise-linear function on the branched cover.  This construction generalizes the usual correspondence between toric Cartier divisors and piecewise-linear functions.  We apply this combinatorial geometric technique to investigate the existence of resolutions of coherent sheaves by vector bundles, using singular nonquasiprojective toric threefolds as a testing ground.  Our main new result is the construction of complete toric threefolds that have no nontrivial toric vector bundles of rank less than or equal to three.  The combinatorial geometric sections of the paper, which develop a theory of cone complexes and their branched covers, can be read independently.
\end{abstract}

\maketitle

\tableofcontents

\section{Introduction}

An algebraic variety $X$ has the resolution property if every coherent sheaf on $X$ is a quotient of a locally free sheaf of finite rank.  If $X$ has the resolution property, then any coherent sheaf $\F$ on $X$ has a possibly infinite resolution by vector bundles
\[
\cdots \rightarrow \E_2 \rightarrow \E_1 \rightarrow \E_0 \rightarrow \F \rightarrow 0.
\]
Totaro asks the following question in much greater generality, for a natural class of stacks that includes all algebraic varieties \cite{Totaro04}, but it remains open and interesting in this special case.

\begin{TQ}
Does every algebraic variety have the resolution property?
\end{TQ}

Roughly speaking, Totaro's question asks whether every algebraic variety has enough nontrivial vector bundles.  The answer is affirmative for quasiprojective varieties, no matter how singular, and for smooth or $\Q$-factorial varieties, no matter how far from quasiprojective \cite{IllusieSGA6}; in both of these cases, every coherent sheaf has a resolution by sums of line bundles.  However, there are many singular nonquasiprojective varieties, starting in dimension 2, with no nontrivial line bundles.

Recently, Schr\"oer and Vezzosi have given an affirmative answer to Totaro's question for normal surfaces \cite{SchroerVezzosi04}.  Given a singular surface $X$, they choose a resolution of singularities $X' \rightarrow X$ and show that the vector bundles on $X'$ that are pulled back from $X$ are exactly those whose restriction to some fixed multiple of the exceptional divisor is trivial.  When $X$ is complete, they use this characterization to construct many nontrivial vector bundles on $X'$ that are pulled back from $X$.  In particular, they show that every complete normal surface has nontrivial vector bundles of rank two.

This paper begins an investigation into Totaro's question for threefolds, using singular nonquasiprojective toric examples.  Our starting point is the observation that, for a complete toric variety $X$, the resolution property implies the existence of nontrivial toric vector bundles, that is, $T$-equivariant vector bundles for the dense torus $T \subset X$ whose underlying vector bundles are nontrivial.  To see this, suppose that $X$ has the resolution property and let $\O_x$ be the structure sheaf of a smooth point $x \in X$, with
\[
\cdots \rightarrow \E_2 \xrightarrow{d_2} \E_1 \xrightarrow{d_1} \E_0 \rightarrow \O_x \rightarrow 0
\]
a resolution by vector bundles.  By the syzygy theorem for regular local rings \cite[Chapter 19]{Eisenbud95}, the kernel of $d_{n-1}$ is locally free, where $n = \dim X$, so we may assume that $\E_i = 0$ for $i > n$.  Then the $K$-theory class
\[
[\E_0] - [\E_1] + \cdots + (-1)^n [\E_n]
\]
is nontrivial; it has nontrivial Chern character by Riemann-Roch for singular varieties \cite[Chapter 18]{IT}.  Since the forgetful map from the Grothendieck group of toric vector bundles to the ordinary Grothendieck group of vector bundles
\[
K^0_T(X) \rightarrow K^0(X)
\]
is surjective \cite[Theorem 6.4]{Merkurjev97}, it follows that $X$ has a toric vector bundle with nontrivial Chern class.

However, there is no known way of constructing a nontrivial toric vector bundle on an arbitrary complete toric variety.

\begin{thm} \label{main}
There are complete toric threefolds with no nontrivial toric vector bundles of rank less than or equal to three.
\end{thm}

\noindent We construct specific such examples in Section \ref{examples}.  It is not known whether these varieties have the resolution property.  It is also not known whether they have any nontrivial vector bundles at all.

\vspace{5 pt}

Although we are still far from resolving Totaro's question, the link between the resolution property and toric vector bundles seems potentially helpful.  The problem of understanding arbitrary vector bundles on algebraic varieties is intractable with current methods, even for vector bundles of small rank on $\P^n$.  For instance, it remains unknown whether every vector bundle of rank 2 on $\P^n$ splits as a sum of line bundles for $n \geq 5$, although Hartshorne's conjecture on complete intersections, which has been open for over thirty years,  would imply that they do split for $n \geq 7$ \cite[p.~89]{OSS}.  The analogous question for toric vector bundles is well-understood; every toric vector bundle of rank less than $n$ on $\P^n$ splits as a sum of line bundles \cite{Kaneyama88}.  Furthermore, Totaro has shown that if every equivariant coherent sheaf on a toric variety $X$ has a resolution by toric vector bundles, then $X$ does have the resolution property \cite[Corollary~5.2]{Totaro04}

\vspace{5 pt}

In order to prove Theorem \ref{main}, we introduce the following combinatorial geometric technique for studying toric vector bundles.  Suppose $X = X(\Delta)$ is the toric variety associated to a fan $\Delta$, over a field $k$.  To each rank $r$ toric vector bundle $\E$ on $X$, we associate a ``branched cover of degree $r$",
\[
\Delta_\E \rightarrow \Delta,
\]
where $\Delta_\E$ is the set of cones of a rational polyhedral cone complex with integral structure, in the sense of \cite[pp.~69--70]{KKMS},  together with an integral piecewise-linear function $\Psi_\E$ on the underlying space $|\Delta_\E|$ of the cone complex.  This $\Delta_\E$ may be described equivalently as a ``multifan", in the sense of \cite{HattoriMasuda03}.  The map $\Delta_\E \rightarrow \Delta$ corresponds to a continuous map
\[
|\Delta_\E| \rightarrow |\Delta|,
\]
taking each cone of $\Delta_\E$ isomorphically onto some cone of $\Delta$, such that the preimage of each cone $\sigma \in \Delta$ is exactly $r$ copies of $\sigma$, counted with multiplicity.  If $r = 1$, then $\E = \O(D)$ for some $T$-Cartier divisor $D$, $\Delta_\E = \Delta$ is the unique branched cover of degree one, and $\Psi_\E = - \Psi_D$ is the negative of the usual piecewise-linear function on $|\Delta|$ associated to $D$ (see \cite[Section 3.4]{Fulton93})\footnote{If $D = \sum d_\rho D_\rho$ is a $T$-Cartier divisor, with $D_\rho$ the $T$-invariant prime divisor corresponding to a ray $\rho$ in $\Delta$, then the usual convention is that $\Psi_D$ evaluated at the primitive generator $v_\rho$ of $\rho$ is equal to $-d_\rho$.  This sign convention is prevalent but not universal, and would lead to a painful proliferation of signs if generalized to toric vector bundles.  For this reason, we follow the minority convention, setting $\Psi_{\O(D)}(v_\rho) = d_\rho$, and defining the line bundle $\L_u$ corresponding to a linear function $u \in M$ to be $\O(\divisor \chi^ {u})$.}.  We illustrate the construction of $\Delta_\E$ and $\Psi_\E$ for a specific rank two example, the tangent bundle to $\P^2$, at the end of the introduction.

The combinatorial data $(\Delta_\E, \Psi_\E)$ are given by equivariant trivializations of $\E$ on the $T$-invariant affine open subvarieties $U_\sigma \subset X$.  This data is not enough to determine the isomorphism class of $\E$ in general; $(\Delta_\E, \Psi_\E)$ determines only the equivariant total Chern class of $\E$, and there may be many nonisomorphic toric vector bundles with the same equivariant total Chern class.   
\vspace{5 pt}

Our strategy for studying toric vector bundles of rank $r$ on $X(\Delta)$ is to divide the problem into two parts.
\begin{enumerate}
\item  Find all of the degree $r$ branched covers of $\Delta$ and compute their groups of piecewise-linear  functions.
\item Study the moduli of rank $r$ toric vector bundles on $X(\Delta)$ with fixed equivariant Chern class.
\end{enumerate}
The key features which make this program workable in relatively simple cases, when $r$ is small and $X(\Delta)$ has few $T$-invariant divisors, are as follows.  First, there are only finitely many degree $r$ branched covers of $\Delta$, and they depend only on the combinatorial type of the fan.  The number of these branched covers may be large, but it is possible to write a computer program to list them all and to compute their groups of piecewise-linear functions.  Second, to understand the moduli, we can consider ``framed" toric vector bundles, that is, toric vector bundles $\E$ with a fixed isomorphism $\E_{x_0} \cong k^r$ of the fiber over the identity $x_0$ in the dense torus with a trivial $r$-dimensional vector space.  There is a fine moduli scheme for framed toric vector bundles with fixed equivariant Chern classes that can be described as a locally closed subscheme of a product of partial flag varieties cut out by rank conditions.  See \cite{moduli} for details.

In the following example, we adopt the na\"{\i}ve viewpoint that a cone complex is a space obtained by gluing a collection of cones along faces, and that a piecewise-linear function on a cone complex is a function that is linear on each cone.  More precise definitions are given in Section \ref{cone complexes}.  We denote cone complexes with bold symbols, and write $|\Delta|$ for the underlying space of a cone complex $\DDelta$.  The example is meant to give an intuition for how we construct a cone complex $\DDelta_\E$, together with a piecewise-linear function $\Psi_\E$ on $|\Delta_\E|$ and a natural projection $|\Delta_\E| \rightarrow |\Delta|$, from a toric vector bundle $\E$ on $X(\Delta)$.  It also illustrates the rough analogy between the topology of this projection of cone complexes and the topology of holomorphic branched covers of Riemann surfaces motivating the terminology ``branched covers of fans".

\begin{example}[Tangent bundle of the projective plane] \label{tangent bundle}
Let $N$ be the rank two lattice
\[ 
N = \Z^3 / (1,1,1).
\]
Let $v_i$ be the image of $e_i$ in $N$, and let $\Delta$ be the fan with maximal cones 
\[
\begin{array}{llll}
\sigma_1 = \<v_2, v_3 \>, & \sigma_2 = \< v_1, v_3 \>, & \mbox{and} & \sigma_3 = \< v_1, v_2 \>.
\end{array}
\]
The associated toric variety $X = X(\Delta)$ is isomorphic to the projective plane.

Let $U_i = U_{\sigma_i}$ be the affine chart corresponding to $\sigma_i$, and let $D_j$ be the prime $T$-invariant divisor on $X$ corresponding to the ray spanned by $v_j$.  On the affine chart $U_1$, the tangent bundle $TU_{1}$ splits $T$-equivariantly as
\[
TU_{1} \cong \O(D_2)|_{U_{1}} \oplus \O(D_3)|_{U_{1}}.
\]
Similarly, $TU_{2} \cong \O(D_1)|_{U_{2}} \oplus \O(D_3)|_{U_{2}}$ and $TU_{3} \cong \O(D_1)|_{U_{3}} \oplus \O(D_2)|_{U_{3}}$.  The linear function associated to $\O(D_j)|_{U_i}$ is $e_j^* - e_i^*$, since $(e_j^* - e_i^*) (v_j) = 1$ and $(e_j^* - e_i^*) (v_k) = 0$, for $i$, $j$, and $k$ distinct.

The cone complex $\DDelta_{TU_{_1}}$ that we associate to $TU_1$ has two maximal cones, $\sigma_{12}$ and $\sigma_{13}$, corresponding to the factors $\O(D_2)|_{U_{_1}}$ and $\O(D_3)|_{U_{1}}$, respectively, and the projection $\pi_1 : |\Delta_{TU_1}| \rightarrow \sigma_1$ maps each of these cones homeomorphically onto $\sigma_1$.  The piecewise-linear function $\Psi_{TU_1}$ is defined by
\[
\Psi_{TU_1}|_{\sigma_{1j}} = \pi_1^* (e_j^* - e_1^*),
\] 
and $\DDelta_{TU_1}$ is defined to the space obtained by identifying the vertices of $\sigma_{12}$ and $\sigma_{13}$, because the vertex is the maximal face of $\sigma_1$ on which the functions $e_2^* - e_1^*$ and $e_3^* - e_1^*$ agree.

\vspace{10 pt}

\begin{center}
\begin{picture}(350,100)(0,0)
\put(30,0){\includegraphics{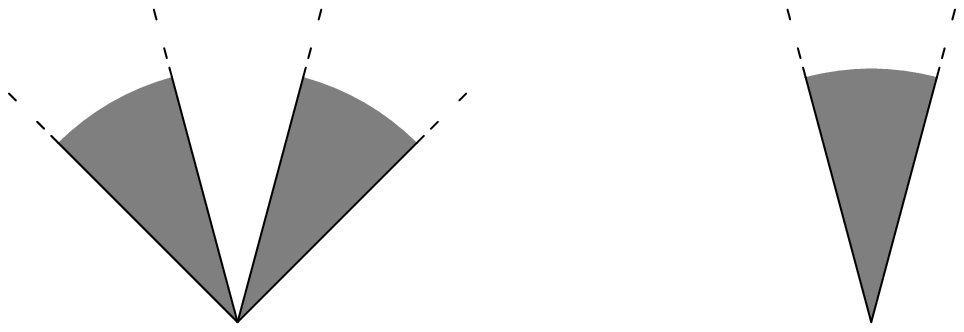}}

\put(30,20){$|\Delta_{TU_1}|$}
\put(310,20){$\sigma_1$}
\put(57,75){$\sigma_{12}$}
\put(147,75){$\sigma_{13}$}
\put(195,40){\vector(1,0){40}}
\put(210,47){$\pi_1$}
\end{picture}
\end{center}

\vspace{10 pt}

Similarly, the cone complex $\DDelta_{TU_2}$ has maximal cones $\sigma_{21}$ and $\sigma_{23}$, and $\DDelta_{TU_3}$ has maximal cones $\sigma_{31}$ and $\sigma_{32}$, and there are projections $\pi_i$ mapping $\sigma_{ij}$ isomorphically onto $\sigma_i$ such that the piecewise-linear functions $\Psi_{TU_i}$ are given by $\Psi_{TU_i}|_{\sigma_{ij}} = \pi_i^*(e_j^* - e_i^*)$.

We construct the cone complex $\DDelta_{TX}$ associated to the tangent bundle $TX$ by gluing together the cone complexes $\DDelta_{TU_i}$ analogously to the way in which $X$ is constructed by gluing the affine charts $U_i$ and the tangent bundle $TX$ is constructed by gluing the $TU_i$.  More precisely, the cone complex $\DDelta_{TX}$ is the cone complex
\[
\DDelta_{TX} = \bigl( \DDelta_{TU_1} \sqcup \DDelta_{TU_2} \sqcup \DDelta_{TU_3} \bigr) / \sim
\]
obtained by gluing $\sigma_{ij}$ to $\sigma_{ji}$ and gluing $\sigma_{ik}$ to $\sigma_{jk}$ along their respective preimages of the ray spanned by $v_k$, for $i$, $j$, and $k$ distinct.  This is the maximal collection of gluings on the $\DDelta_{TU_i}$ that is compatible with the piecewise-linear functions $\Psi_{TU_i}$ and the projections to $N_\R$.  The piecewise-linear function $\Psi_{TX}$ is defined to be the function whose restriction to $|\Delta_{TU_i}|$ is $\Psi_{TU_i}$.

In summary, $\DDelta_{TX}$ has six maximal cones $\sigma_{ij}$ for $1 \leq i,j \leq 3$ and $i \neq j$, the projection $\pi: |\Delta_{TX}| \rightarrow |\Delta|$ maps $\sigma_{ij}$ homeomorphically onto $\sigma_i$, and the piecewise-linear function $\Psi_{TX}$ is given by
\[
\Psi_{TX}|_{\sigma_{ij}} = \pi^*(e_j^* - e_i^*).
\]
The projection may be drawn as follows, with $\sigma_{ij}$ shaded light for $i-j \equiv 1$ mod 3, and dark for $i-j \equiv 2$ mod 3.

\vspace{10 pt}

\begin{center}
\begin{picture}(350,100)(0,0)

\put(40,0){\includegraphics{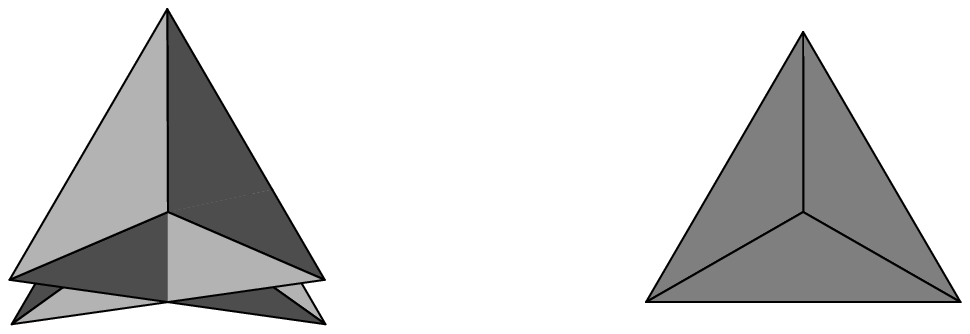}}

\put(165,40){\vector(1,0){40}}
\put(180,47){$\pi$}

\put(10,50){$|\Delta_{TX}|$}
\put(320,50){$|\Delta|$}

\end{picture}
\end{center}

\vspace{10 pt}

\noindent It is not difficult to see that, as topological spaces, $|\Delta_{TX}|$ and $|\Delta|$ are homeomorphic to the complex plane, and the homeomorphisms may be chosen such that the projection $\pi: |\Delta_{TX}| \rightarrow |\Delta|$ is identified with the branched cover of Riemann surfaces $z \mapsto z^2$.

\end{example}

\vspace{10 pt}

\noindent \textbf{Acknowledgments.} I have benefited from many conversations with friends and colleagues while preparing this work.  I thank, in particular, D.~Abramovich, A.~Buch,  I.~Dolgachev, W.~Fulton, J.~Harris, M.~Hering, R.~Lazarsfeld, M.~Musta\c{t}\v{a}, B.~Totaro, and R.~Vakil for helpful comments and suggestions.  

I wrote this paper while at the University of Michigan, Ann Arbor and am grateful to the Department of Mathematics there for providing a wonderful research environment.

\section{Cone complexes and branched covers} \label{cone complexes}

Cone complexes are an intuitive and flexible framework for the application of elementary geometric, topological, and combinatorial techniques; they appear naturally in a variety of contexts related to the geometry of real and complex torus actions.  The definition that we use was given by Kempf, Knudsen, Mumford, and Saint-Donat over thirty years ago in connection with their classification of singularities of toroidal embeddings \cite{KKMS}. Equivalent forms of this definition continue to be used in the study of logarithmic schemes; see, for instance, \cite{Kato94} and \cite{AbramovichKaru00}.  Instances of cone complexes have also appeared elsewhere in the literature, implicitly, under various names: ``fans'' in the classification of toric varieties, ``multifans'' in the classification of torus manifolds and torus orbifolds \cite{HattoriMasuda03}, ``systems of fans'' in the classification of toric prevarieties \cite{ACampoNeuenHausen01}, and now ``branched covers of fans'' in the classification of toric vector bundles.

In Sections~\ref{definitions section} and \ref{morphisms section}, we systematically develop a basic categorical theory of cone complexes and their morphisms, since this is necessary for defining branched covers, and there seems to be no suitable reference in the literature.  We then introduce branched covers of cone complexes from topological and combinatorial points of view in Sections~\ref{top section} and \ref{poset section}, respectively.  The essential equivalence between these viewpoints is established in Proposition~\ref{top to poset}.  The degree $d$ branched covers of a fixed cone complex carry a natural partial order, given in Definition~\ref{branched order}.  The main result that we will use in the analysis of toric vector bundles on specific toric threefolds in Section~\ref{examples} is a bijective correspondence between branched covers of a complete 3-dimensional fan that are maximal in this partial order and holomorphic branched covers of $\P^1(\C)$ with specified branch locus.  See Theorem~\ref{bijections}.

\vspace{5 pt}

\subsection{Rational polyhedral cone complexes with integral structure} \label{definitions section}

An integral vector space is a real vector space $N_{\R}$, together with a lattice $N \subset N_\R$ such that $N_{\R} = N \otimes_{\Z} \R$.  We write $M$ for the dual lattice $M = N^{\vee}$.

Roughly speaking, a rational polyhedral cone complex with integral structure is a collection of rational convex polyhedral cones in integral vector spaces glued together along isomorphic faces.  For brevity, we write ``cone complex'' to mean ``rational polyhedral cone complex with integral structure.''\footnote{All polyhedral cone complexes considered in this paper are rational and have an integral structure.  Following \cite{KKMS}, one can also define rational polyhedral cone complexes without integral structure as triples $( |\Delta|, \Delta, \{ M_{\sigma} \}_{\sigma \in \Delta} )$, where each $M_{\sigma}$ is a finite dimensional rational vector space of continuous functions on $\sigma$.  Other natural variations are possible; for instance, one could allow the $\phi_{\sigma}(\sigma)$ to be not-necessarily-rational polyhedral cones, and let the $M_{\sigma}$ be abelian groups, rational vector spaces, or real vector spaces of continuous functions.} The following definition is adapted from \cite[pp.~69--70]{KKMS}.\footnote{There are two main differences, essentially cosmetic, between Definition \ref{cone complex} and the definition given by Mumford and his collaborators.  First, we explicitly require that $\phi_\sigma(\sigma)$ be a rational cone, which seems to have been implicit in the original definition.  Second, we omit the requirement that $\phi_\sigma(\sigma)$ be full-dimensional, since this follows from the hypothesis that $M_\sigma$ is a subgroup of the continuous functions on $\sigma$.}

\begin{defn} \label{cone complex}
A cone complex $\DDelta$ is a topological space $|\Delta|$ together with a finite collection $\Delta$ of closed subspaces of $|\Delta|$ and, for each $\sigma \in \Delta$, a finitely generated group $M_{\sigma}$ of continuous functions on $\sigma$, satisfying the following conditions.  Let $N_{\sigma} = M_{\sigma}^{\vee}$.
\begin{enumerate}
\item The natural map $\phi_{\sigma}: \sigma \rightarrow (N_{\sigma})_{\R}$ given by $x \mapsto [u \mapsto u(x)]$ maps $\sigma$ homeomorphically onto a rational convex polyhedral cone.
	
\item The preimage under $\phi_{\sigma}$ of each face of $\phi_{\sigma}(\sigma)$ is some $\tau \in \Delta$, and $M_{\tau} = \{ u|_{\tau} : u \in M_{\sigma} \}$.

\item The topological space $|\Delta|$ is the disjoint union of the relative interiors of the $\sigma \in \Delta$.
\end{enumerate}
\end{defn}

We say that $\Delta$ is the set of cones of $\DDelta$.  If $\tau$ is the preimage in $|\Delta|$ of a face of $\phi_{\sigma}(\sigma)$ then we say that $\tau$ is a face of $\sigma$ and write $\tau \preceq \sigma$.  We write $\tau \prec \sigma$ if $\tau$ is a face of $\sigma$ and $\tau \neq \sigma$.  We say that $\sigma \in \Delta$ is strongly convex (resp.\ simplicial, $i$-dimensional) if $\phi_{\sigma}(\sigma)$ is, and that $\DDelta$ is strongly convex (resp.\  simplicial) if all of its cones are.

\begin{example}[Integral vector spaces, cones, and fans]
Suppose $N_{\R}$ is an integral vector space.  Then $\mathbf{N_{\R}} = (N_{\R}, \{ N_{\R} \}, N^{\vee} )$ is a cone complex.  

For each convex rational polyhedral cone $\sigma \subset N_{\R}$, define $M_{\sigma} = M/(M \cap \sigma^{\perp})$.  Then $\pmb{\sigma} = (\sigma, \{ \tau \preceq \sigma \}, \{ M_{\tau}\})$ is a cone complex.  Furthermore, if $\Delta$ is a fan in $N_{\R}$, then $\DDelta = (|\Delta|, \Delta, \{M_{\sigma}\}_{\sigma \in \Delta})$ is a strongly convex cone complex.
\end{example}

\begin{example}[Puckered cone] \label{puckered cone}
Let $\sigma$ be the positive quadrant in $\R^2$.  The cone complex obtained by gluing two copies of $\sigma$ along both axes, which we call a ``puckered cone", is the multifan corresponding to the ``bug-eyed plane", the nonseparated toric prevariety consisting of two copies of $\A^2$ glued along $\A^2 \smallsetminus 0$.  See \cite[Section 4]{chowcohom} for a discussion of multifans and toric prevarieties.
\end{example}

\noindent A map from the puckered cone in Example \ref{puckered cone} to $\R^n$ that is linear on each cone cannot be injective.  Such nonembeddable cone complexes occur naturally in many situations, including in the proofs of Lemma \ref{max cones} and Proposition \ref{codim 1}, below.

\begin{example}[Cone over a simplicial complex] \label{simplicial complex}
Let $X$ be the geometric realization of a simplicial complex.  Then $\CX$, the cone over $X$, is the cone complex whose underlying topological space is
\[
|CX| = X \times [0,\infty) / (X \times \{ 0 \} ),
\]
and whose cones are the cones over the simplices in $X$, where the integral linear functions on the cone over a simplex are the linear functions taking integer values at $(v,1)$ for each vertex $v$ of the simplex. 
\end{example}

\noindent Many more simple examples of cone complexes may be constructed from the examples above by sequences of disjoint unions, products, and subdivisions.  

Another important example, which can be helpful for thinking about branched covers (see Remark \ref{branched stars}), is the star of a cone in a cone complex.  Roughly speaking, the star of a cone is what you see looking out from a point in its relative interior.  For a cone $\sigma \succeq \tau$, let $\sigma_\tau$ be the rational polyhedral cone in $(N_\sigma)_\R$
\[
\sigma_\tau = \sigma + (-\tau),
\]
which is the cone spanned by $\sigma - v$ for any point $v$ in the relative interior of $\tau$.

\begin{example}[Star of a cone]
Let $\DDelta$ be a cone complex, with $\tau \in \Delta$ a cone.  Then
\[
\Star(\tau) = \Bigl( \, \bigsqcup_{\sigma \succeq \tau} \sigma_\tau \Bigr) \,  \Big/ \sim \, ,
\]
the space obtained identifying $\gamma_\tau$ with its image in $\sigma_\tau$, for $\tau \preceq \gamma \preceq \sigma$, is a cone complex.
\end{example}

\begin{remk} \label{connected components}
A cone in $\Delta$ is minimal if and only if it has no proper faces, so the minimal cones are exactly those $L \in \Delta$ such that $\phi_{L}(L) = (N_{L})_{\R}$.  If $L$ is a minimal cone, then $\Star(L)$ is naturally identified with the subcomplex of $\DDelta$ that is the union of the cones containing $L$. Since every cone in $\Delta$ contains a unique minimal cone and the star of each minimal cone is closed and connected, the connected components of a cone complex $\DDelta$ are exactly those $\Star(L)$ such that $L$ is a minimal cone in $\Delta$.  Therefore, a cone complex is connected if and only if it contains a unique minimal cone $L_\Delta$.  In particular, a strongly convex, connected cone complex $\DDelta$ contains a unique zero-dimensional cone $0_\Delta$.
\end{remk}

\begin{example}[Wedge sums]
Suppose $\DDelta$ and $\DDelta'$ are strongly convex, connected cone complexes.  The wedge sum $\DDelta \vee \DDelta'$ is the cone complex obtained from the disjoint union $\DDelta \sqcup \DDelta'$ by identifying $0_{\Delta}$ and $0_{\Delta'}$.
\end{example}

\begin{defn}[Piecewise-linear functions]
Let $\DDelta$ be a cone complex.  Then the group of integral piecewise-linear functions on $\DDelta$ is
\[
\PL(\DDelta) = \{ \Psi: |\Delta| \rightarrow \R \ | \ \Psi|_\sigma \in M_\sigma \mbox{ for all } \sigma \in \Delta \}.
\]
\end{defn}

\noindent If $\Delta$ is a fan, then $\PL(\DDelta)$ is naturally identified with the group of $T$-Cartier divisors on the toric variety $X(\Delta)$ \cite[Section~3.4]{Fulton93}.

\vspace{0 pt}

\subsection{Morphisms of cone complexes} \label{morphisms section}

Roughly speaking, a morphism of cone complexes should be a continuous, piecewise-linear map that is linear on each cone and respects the integral
structures.  We define a morphism of cone complexes as follows.

\begin{defn}
A morphism of cone complexes $\bf: \DDelta' \rightarrow \DDelta$ is a continuous map of topological spaces $|\Delta'| \rightarrow |\Delta|$ such that, for each cone $\sigma' \in \Delta'$, there is some $\sigma \in \Delta$ with $\bf(\sigma') \subset \sigma$ and $\bf^{*} M_{\sigma} \subset M_{\sigma'}$.
\end{defn}

\noindent The composition of any two morphisms of cone complexes is a morphism, and the identity is a morphism from any cone complex to itself; cone complexes with morphisms of cone complexes form a category.  

If $\bf: \DDelta' \rightarrow \DDelta$ is a morphism, then $\bf^*(\PL(\DDelta)) \subset \PL(\DDelta')$, so $\PL$ is a contravariant functor from cone complexes to finitely generated free abelian groups.

\begin{example}
Suppose $\Delta$ is a fan in $N_{\R}$.  The natural inclusion $\DDelta \hookrightarrow \mathbf{N}_{\R}$ is a morphism of cone complexes.  If $\Delta'$ is a fan in $N'_{\R}$ and $\varphi: N' \rightarrow N$ induces a morphism of fans, then $\varphi_{\R}|_{\DDelta'} : \DDelta' \rightarrow \DDelta$ is a morphism of cone complexes.
\end{example}

\begin{example} In Example \ref{puckered cone}, the projection from the puckered cone to $\sigma$ is a morphism of cone complexes.
\end{example}

\begin{example}[Fibered products]
If $\DDelta$ and $\DDelta'$ are cone complexes with morphisms to a cone complex $\SSigma$, then the fibered product $\DDelta \times_{\SSigma} \DDelta'$ exists as a cone complex whose underlying topological space is $|\Delta| \times_{|\Sigma|} |\Delta'|$.  The cones are the fibered products
\[
\{ \sigma \times_{|\Sigma|} \sigma' \ | \ \sigma \in \Delta \mbox{ and } \sigma' \in \Delta' \},
\]
with integral linear functions given by 
\[
M_{\sigma \times_{|\Sigma|} \sigma'} = \{ p_{1}^{*} (u) + p_{2}^{*} (u') \ | \ u \in M_{\sigma} \mbox{ and } u' \in M_{\sigma'} \}.
\]
The projections to each factor are morphisms of cone complexes over $\SSigma$, with the usual universal property.
\end{example}      

\begin{remk}
If $\Delta$ and $\Delta'$ are fans in $N_{\R}$ with the same support, then $\DDelta \times_{\mathbf{N_{\R}}} \DDelta'$ is the smallest common refinement of $\DDelta$ and $\DDelta'$, and the projections correspond to the proper birational toric morphisms induced by the identity on $N$.
\end{remk}

If $\bf: \DDelta \rightarrow \SSigma$ is a morphism of cone complexes, then $\bf$ maps the relative interior of a cone $\sigma \in \Delta$ into the relative interior of a unique cone $f(\sigma) \in \Sigma$.

\begin{example}
Let $\bf: \DDelta \rightarrow \SSigma$ be a morphism.  Then $\bf|_\sigma : \sigma \rightarrow f(\sigma)$ is the restriction of a linear map $(N_\sigma)_\R \rightarrow (N_{f(\sigma)})_\R$ that takes $\sigma_\tau$ into $f(\sigma)_{f(\tau)}$ for each face $\tau$ of $\sigma$.  For fixed $\tau \in \Delta$, these maps glue together to give a morphism
\[
\Star_\tau (\bf) : \Star(\tau) \rightarrow \Star(f(\tau)).
\]
\end{example}

\vspace{0 pt}

\subsection{Branched covers of cone complexes} \label{top section}

The theory of branched covers of cone complexes that we introduce here is loosely based on the classical theory of branched covers of Riemann surfaces.  We begin our study with a motivational example.  Recall that if $f : Y \rightarrow X$ is a map of sets such that $f^{-1}(x)$ is finite for each $x \in X$ then, for a function $w$ on $Y$, we define the trace of $w$ to be the function on $X$ defined by
\[
   \Tr_{f} (w) (x) = \sum_{f(y) = x} w(y).
\]

\begin{example}[Triangulated branched covers of Riemann surfaces] \label{triangulated surfaces}
Let $f:Y \rightarrow X$ be a degree $d$ branched cover of compact Riemann surfaces.  In other words, $X$ and $Y$ are compact two-dimensional manifolds, each with a given complex analytic structure, and for each point $y$ in $Y$ there is a positive integer $w(y)$ such that $f$ is given in local coordinates around $y$ by
\[
z \mapsto z^{w(y)},
\]
and such that $\Tr_{f}(w) = d$.  Because $Y$ is compact, $w$ is equal to one at all but finitely many points of $Y$.  The points where $w$ is greater than one are called ramification points, and the images of the ramification points in $X$ are called branch points.

Fix a triangulation of $X$ such that every branch point is a vertex and no two branch points are joined by an edge.  This triangulation lifts to a triangulation of $Y$, so we may think of $X$ and $Y$ as simplicial complexes, with $f$ a simplicial map.  Let $\CX$ and $\CY$ be the cones over the simplicial complexes $X$ and $Y$, respectively, as in example \ref{simplicial complex}.  The function $w$ on $Y$ induces a function $\bw$ on $\CY$ by setting $\bw((y,t)) = w(y)$ for $t \neq 0$, and $\bw(0_{\CY}) = d$.  Note that $\bw$ is constant on the relative interior of each cone of $CY$.  The simplicial map $f$ induces a morphism of cone complexes $\bf: \CY \rightarrow \CX$, given by $(y,t) \mapsto (f(y),t)$, that has the following two properties. 
\begin{enumerate}

\item For each cone $\sigma$ in $CY$, $\bf$ maps $\sigma$ isomorphically onto $f(\sigma)$.

\item For each connected open set $U \subset |CX|$, and for each connected component $V \subset \bf^{-1}(U)$, $\Tr_{\bf|_V}(\bw)$ is a constant function on $U$.
\end{enumerate}
\end{example}

We extrapolate from Example \ref{triangulated surfaces} to define branched covers of cone complexes as follows.

\begin{defn}
A weighted cone complex is a cone complex $\DDelta$ together with a function $\bw : |\Delta| \rightarrow \Z^{+}$ that is constant on the relative interior of each cone.
\end{defn}

\begin{defn}[Branched covers] \label{top covers}
Let $\SSigma$ be a connected cone complex, and let $(\DDelta, \bw)$ be a weighted connected cone complex.  A branched cover $\bf: \DDelta \rightarrow \SSigma$ is a morphism of cone complexes such that
\begin{enumerate}

\item For each cone $\sigma$ in $\Delta$, $\bf$ maps $\sigma$ isomorphically onto $f(\sigma)$.

\item For each connected open set $U \subset |\Sigma|$, and for each connected component $V \subset \bf^{-1}(U)$, $\Tr_{\bf|_V}(\bw)$ is a constant function on $U$.
\end{enumerate}
\end{defn}

\noindent If $\bf: \DDelta \rightarrow \SSigma$ is a branched cover, then $\Tr_\bf(\bw)$ is a positive integer, which we call the degree of $\bf$.  The degree of $\bf$ is equal to $\bw(v)$ for any $v$ in the relative interior of the minimal cone $L_\Delta$.  We say that  $\bf$ is ramified along those $\tau \neq L_\Delta$ such that $\bw > 1$ on the relative interior of $\tau$.

\begin{example}
Let $|\Delta_{T\P^2}|$ and $|\Delta|$ be as in Example \ref{tangent bundle}, and let $\DDelta_{T\P^2}$ and $\DDelta$ be the natural cone complexes with underlying topological spaces $|\Delta_{T\P^2}|$ and $|\Delta|$, respectively.  Let $\DDelta_{T\P^2}$ be weighted by $\bw$, with $\bw = 1$ everywhere except at the vertex, where $\bw(0_{\Delta_{T\P^2}}) = 2$.  Then the projection $\pi: \DDelta_{T\P^2} \rightarrow \DDelta$ is a branched cover of degree 2.  
\end{example}

\begin{example}[$d$-fold puckered cone] \label{d-fold pucker}
Let $\sigma$ be a cone, and let $\Puck^d(\sigma)$ be the cone complex obtained by gluing $d$ copies of $\sigma$ along their entire boundaries.  Let $\bw$ be the weight function on $\Puck^d(\sigma)$ that is 1 on the interior of each copy of $\sigma$, and $d$ along their shared boundary.   Then the natural projection $\pi: \Puck^d(\sigma) \rightarrow \pmb{\sigma}$ is a branched cover of degree $d$.  
\end{example}

\noindent Example \ref{d-fold pucker} generalizes the puckered cone in Example \ref{puckered cone}, which is the case where $d = 2$ and $\sigma$ is the positive quadrant in $\R^2$.

The following are a few other simple constructions of branched covers.  If $d$ is a positive integer and $d \DDelta$ is the cone complex $\DDelta$ weighted by the constant function $\bw = d$, then the identity on $\DDelta$ induces a degree $d$ branched cover $d \DDelta \rightarrow \DDelta$.

If $\DDelta$ is strongly convex, let $\vee^d \DDelta$ be the wedge sum of $d$ copies of $\DDelta$, weighted by $\bw$, with $\bw = 1$ everywhere except at the vertex, where $\bw(0_{\vee^d \DDelta}) = d$.  Then the natural projection $\vee^d \DDelta \rightarrow \DDelta$ is a degree $d$ branched cover.  

More generally, the wedge sum of any two branched covers of the same cone complex, of degree $d$ and $d'$, respectively, is a branched cover of degree $d + d'$.  Similarly, if $(\DDelta,\bw)$ and $(\DDelta',\bw')$ are branched covers of $\SSigma$ of degrees $d$ and $d'$, respectively, then let $\bw \cdot \bw'$ be the weight function on $\DDelta \times_{\SSigma} \DDelta'$ given by 
\[
\bw \cdot \bw' ((x,y)) = \bw(x) \cdot \bw'(y).
\]
Then $(\DDelta \times_{\SSigma} \DDelta', \bw \cdot \bw')$ is a branched cover of $\SSigma$ of degree $d \cdot d'$. 

\begin{remk}
In \cite{HattoriMasuda03}, a weighted multifan in $N_\R$ is called complete if the projection to $N_\R$ satisfies the trace condition in Definition \ref{top covers}.  In particular, a branched cover of a complete fan is a complete multifan.
\end{remk}

\subsection{Branched covers of posets} \label{poset section}

Here we introduce a combinatorial notion of branched covers of posets and give a bijective correspondence between branched covers of a connected cone complex $\SSigma$ and branched covers of its poset of cones.  The rough idea behind this correspondence is that a branched cover $(\DDelta, \bw) \rightarrow \SSigma$ is given by purely combinatorial data: a collection of isomorphic copies of cones in $\Sigma$ with some data for gluing these cones along faces, both encoded in the map of posets $\Delta \rightarrow \Sigma$, plus a weight function.  The resulting bijection is useful for applications, since the cone complexes give geometric intuition and permit the application of topological techniques, as in Section~\ref{topology of maximal}, while the posets are often easier to work with in combinatorial situations.  For instance, in Section~\ref{branched covers and chern classes} we will define the branched cover associated to a toric vector bundle in terms of a branched cover of posets.  Furthermore, we get a natural bijection between branched covers of combinatorially equivalent fans that streamlines many numerical computations.  See, for instance, Examples~\ref{Fulton's threefold} and \ref{new example}.

For a finite partially ordered set $\Po$, with $x \in \Po$, let 
\[
\Po_x = \{ y \in \Po \ | \ y \preceq x \},
\]
and give $\Po$ the poset topology, where $S \subset \Po$ is closed if and only if $S$ contains $\Po_x$ for every $x \in S$.  If $\DDelta$ is a cone complex, then the poset topology on $\Delta$ is the quotient topology for the projection $|\Delta| \rightarrow \Delta$ taking a point $v \in |\Delta|$ to the unique cone containing $v$ in its relative interior.

If $\bf : \DDelta \rightarrow \SSigma$ is a morphism of cone complexes, then the induced map
\[
f: \Delta \rightarrow \Sigma,
\]
taking a cone $\sigma \in \Delta$ to the smallest cone $f(\sigma) \in \Sigma$ containing its image under $\bf$, is an order-preserving, and hence continuous, map of posets.

Let $\SSigma$ be a cone complex, and let $f: \Po \rightarrow \Sigma$ be a map of posets such that $f$ maps $\Po_x$ isomorphically onto $\Sigma_{f(x)}$ for every $x \in \Po$.  Then the fibered product
\[
|\Sigma| \times_\Sigma \Po
\]
is the underlying space of a cone complex $\SSigma \times_\Sigma \Po$ whose cones are
\[
\{ |\Sigma| \times_\Sigma \Po_x \ | \ x \in \Po \}.
\] 
In particular, the poset of cones of $\SSigma \times_\Sigma \Po$ is canonically identified with $\Po$.  The projection $p_1: \SSigma \times_\Sigma \Po \rightarrow \SSigma$ restricts to an isomorphism on each cone of $\SSigma \times_\Sigma \Po$;  it maps $|\Sigma| \times_\Sigma \Po_x$ isomorphically onto $f(x)$ for every $x \in \Po$.

\begin{defn}
A weighted poset is a poset $\Po$ together with a function
\[
w: \Po \rightarrow \Z^+.
\]
\end{defn}

\noindent If $\bw$ is a weight function on a cone complex $\DDelta$, then $\bw$ descends to a weight function $w$ on $\Delta$, by setting $w(\sigma) = \bw(v)$ for any point $v$ in the relative interior of $\sigma$.  Conversely, any weight function on $\DDelta$ is the pullback of a unique function on $\Delta$.

\begin{remk} \label{poset opens}
For a finite poset $\Po$, with $x \in \Po$, note that
\[
\Po^x = \{ y \in \Po \ | \ y \succeq x \}
\]
is the smallest open neighborhood of $x$ in the poset topology.  Furthermore, if $f: \Po \rightarrow \Sigma$ is  a map of posets mapping $\Po_x$, isomorphically onto $\Sigma_{f(x)}$ for all $x \in \Po$, then the connected components of $f^{-1}(\Sigma^\tau)$ are exactly those $\Po^x$ for $x \in f^{-1}(\tau)$.
\end{remk}

Recall that a poset is rooted if it has a unique minimal element.  By Remark~\ref{connected components}, a cone complex $\DDelta$ is connected if and only if its poset of cones $\Delta$ is rooted.

\begin{defn} \label{branched posets}
A branched cover of a rooted poset $\Sigma$ is a weighted rooted poset $(\Po, w)$ with a map $f: \Po \rightarrow \Sigma$ such that, for every $x$ in $\Po$,
\begin{enumerate}
\item $f$ maps $\Po_x$ isomorphically onto $\Po_{f(x)}$.

\item $\Tr_{f|_{\Po^x}}(w)$ is a constant function on $\Sigma^{f(x)}$.
\end{enumerate}
\end{defn}

The essential equivalence between branched covers of posets and branched covers of cone complexes may be expressed as follows.

\begin{prop}  \label{top to poset}
Let $\SSigma$ be a connected cone complex.
\begin{enumerate}

\item If $\bf: (\DDelta, \bw) \rightarrow \SSigma$ is a branched cover of cone complexes then the induced map
\[
f : (\Delta, w) \rightarrow \Sigma
\]
is a branched cover of posets.

\item If $(\Po, w) \rightarrow \Sigma$ is a branched cover of posets then the projection
\[
p_1 : \bigl( \SSigma \times_\Sigma \Po, \, p_2^* w \bigr) \rightarrow \SSigma
\]
is a branched cover of cone complexes.
\end{enumerate}
\end{prop}

\begin{proof}
We prove part (1) of the proposition; the proof of part (2) is similar.  Since $\bf$ is a branched cover, any cone $\tau$ in $\Delta$ maps isomorphically onto $f(\tau)$, so $\Delta_\tau$, the poset of faces of $\tau$, maps isomorphically onto $\Sigma_{f(\tau)}$, the poset of faces of $f(\tau)$.  Now, let $V^\tau \subset |\Delta|$ be the connected open set
\[
V^\tau = \bigcup_{\sigma \succeq \tau} \relint \sigma.
\]
Then the constant function $\Tr_{f|_{V^\tau}}(\bw)$ is the pullback of $\Tr_{f|_{\Delta^\tau}}(w)$ under the projection $V^\tau \rightarrow \Delta^\tau$.  Therefore, $\Tr_{f|_{\Delta^\tau}}(w)$ is a constant function on $\Sigma^{f(\tau)}$, as required.
\end{proof}

Roughly speaking, Proposition \ref{top to poset} gives a correspondence between branched covers of a connected cone complex $\SSigma$ and branched covers of its poset of cones $\Sigma$, up to isomorphisms.  

\begin{defn}
Two branched covers $\bf: (\DDelta, \bw) \rightarrow \SSigma$ and $\bf': (\Delta', \bw') \rightarrow \SSigma$ are isomorphic if there is an isomorphism of cone complexes
\[
\bg: \DDelta' \xrightarrow{\sim} \DDelta
\] 
such that $\bf' = \bf \circ \bg$ and $\bg^* \bw = \bw'$.  
\end{defn}

\noindent Similarly, we say that two poset branched covers $f: (\Po,w) \rightarrow \Sigma$ and $f':(\Po', w') \rightarrow \Sigma$ are isomorphic if there is an isomorphism of posets $g : \Po' \xrightarrow{\sim} \Po$ such that $f' = f \circ g$ and $g^* w = w'$.

The set of isomorphism classes of degree $d$ branched covers of $\SSigma$ is partially ordered as follows. 

\begin{defn} \label{branched order}
Let $\bf:(\DDelta, \bw) \rightarrow \SSigma$ and $\bf' : (\DDelta', \bw') \rightarrow \SSigma$ be branched covers of degree $d$.  Then $\bf' \geq \bf$ if there is a map of cone complexes $\bg : \DDelta' \rightarrow \DDelta$ such that $\bf' = \bf \circ \bg$ and $\Tr_\bg(\bw') = \bw$.  
\end{defn}

\noindent If $\bf' \geq \bf$ then we think of $\DDelta'$ as being constructed by ``pulling apart" some of the ramified cones in $\DDelta$.  For instance,  the $d$-fold puckered cone $\P^d(\sigma)$ in Example \ref{d-fold pucker} may be thought of as being constructed by pulling apart the layers in $d \pmb{\sigma}$.

Similarly, the set of isomorphism classes of degree $d$ poset branched covers is partially ordered by saying that $f' \geq f$ if there is a map of cone complexes $g: \Po' \rightarrow \Po$ such that $f' = f \circ g$ and $\Tr_g(w') = w$.

\begin{prop} \label{top poset equivalence}
Let $\SSigma$ be a connected cone complex and let $d$ be a positive integer.  There is an order-preserving bijection between isomorphism classes of degree $d$ branched covers of $\SSigma$ and isomorphism classes of degree $d$ poset branched covers of $\Sigma$.
\end{prop}

\begin{proof}
The maps given by parts (1) and (2) or Proposition~\ref{top to poset} preserve isomorphism classes, respect the partial orderings, and are inverse to each other on isomorphism classes.
\end{proof}

\begin{remk} \label{branched stars}
It is sometimes helpful to think of a morphism of connected cone complexes $\bf: (\DDelta, \bw) \rightarrow \SSigma$ ``locally", in terms of the maps $\Star_\tau(\bf)$ for $\tau \in \Delta$ and $\tau \neq 0_\Delta$.  The cone complex $\Star(\tau)$ has a natural weight function $\bw_\tau$, given by
\[
w_\tau(\sigma_\tau) = w(\sigma).
\]
It is straightforward to check that a morphism of connected cone complexes $\bf$ is a branched cover of degree $d$ if and only if $\Star_\tau(\bf)$ is a branched cover for every $\tau \neq 0_\Delta$, and $\Tr_f(w) = d$.
\end{remk}

\vspace{5 pt}

\subsection{Topology of maximal branched covers} \label{topology of maximal}

If $\Sigma$ is a fan then, for the purposes of studying rank $r$ toric vector bundles on $X(\Sigma)$, we want to understand the group of piecewise-linear functions $\PL(\DDelta)$ for every degree $r$ branched cover $(\DDelta, \bw) \rightarrow \SSigma$.  Since piecewise-linear functions pull back under morphisms of cone complexes, it is generally enough to understand $\PL(\DDelta)$ for branched covers that are maximal in the partial order given in Definition \ref{branched order}.

Recall that a branched cover $(\DDelta, \bw) \rightarrow \SSigma$ is said to be ramified along a cone $\tau \in \Delta$ if $w(\tau) > 1$ and $\tau$ is not  the minimal cone $L_\Delta$.

\begin{lem} \label{max cones}
Let $\bf: (\DDelta,\bw) \rightarrow \SSigma$ be a maximal branched cover of cone complexes.  Then $\bf$ is unramified along every maximal cone in $\Delta$.
\end{lem}

\begin{proof}
Suppose $\sigma \neq L_\Delta$ is a maximal cone in $\Delta$ with $w(\sigma) = d$.  We must show that $d = 1$.  Let $\DDelta'$ be the cone complex constructed from $\DDelta$ by removing $\sigma$ and gluing in a copy of the $d$-fold puckered cone $\Puck^d(\sigma)$ (see Example \ref{d-fold pucker}), with $\bg: \DDelta' \rightarrow \DDelta$ the natural projection.  Then $\bf' = \bf \circ \bg$ is also a degree $d$ branched cover, and $\bf' \geq \bf$.  Since $\bf$ is maximal, $\bg$ is an isomorphism, and therefore $d = 1$.
\end{proof}

\begin{prop} \label{codim 1}
Let $\Sigma$ be a complete fan, with $\bf: (\DDelta, \bw) \rightarrow \SSigma$ a maximal branched cover.  Then $\bf$ is unramified along every codimension one cone.
\end{prop}

\begin{proof}
Let $\tau \neq 0_\Delta$ be a codimension one cone, with $w(\tau) = d$.  We must show that $d = 1$.  Since $\Sigma$ is a complete fan, $f(\tau)$ is contained in exactly two maximal cones of $\Sigma$, say $\gamma$ and $\gamma'$.  By Lemma \ref{max cones}, there must be exactly $d$ cones $\sigma_1, \ldots, \sigma_d$ (resp.\ $\sigma'_1, \ldots, \sigma'_d$) containing $\tau$ such that $f(\sigma_i) = \gamma$ (resp. $f(\sigma'_i) = \gamma'$).

Let $\Delta'$ be the poset
\[
\Delta' = \bigl( \Delta \smallsetminus \{ \tau \} \bigr) \sqcup \{ \tau_1, \ldots, \tau_d \},
\]
with the following partial ordering.  The $\tau_i$ are pairwise incomparable, 
\begin{eqnarray*}
\tau_i \prec \sigma_j & \mbox{ if and only if } & i = j, \\
\tau _i \prec \sigma'_j & \mbox{ if and only if } & i = j, \\
\rho \prec \tau_i \mbox{ in } \Delta' & \mbox{ if and only if } & \rho \prec \tau \mbox{ in } \Delta,
\end{eqnarray*}
 and all other relations are exactly as in $\Delta$.  
 
Let $w'$ be the weight function on $\Delta'$ given by $w'|_{\Delta \smallsetminus \{ \tau \} } = w$ and $w'(\tau_i) = 1$, with $g: \Delta' \rightarrow \Delta$ the projection taking $\tau_i$ to $\tau$.  Then $f' = f \circ g$ is a branched cover of degree $d$, and $f' \geq f$.  Since $f$ is maximal, $g$ is an isomorphism, and therefore $d = 1$.
\end{proof}

\begin{example}
Let $\Sigma$ be the unique complete one-dimensional fan, and let $\rho$ and $\rho'$ be the rays of $\Sigma$.  A degree $d$ branched cover of $\Sigma$ corresponds to an ordered pair $(\lambda, \lambda')$ of partitions of $d$, such that the weights of the preimages of $\rho$ and $\rho'$ in $\Sigma_{(\lambda, \lambda')}$ are given by $\lambda$ and $\lambda'$, respectively.  The unique maximal degree $d$ branched cover of $\Sigma$, up to isomorphism, is $\vee^d \Sigma$, which corresponds to the pair of partitions $(\lambda, \lambda') = \bigl(1^d, 1^d \bigr)$.
\end{example}

\begin{example}
Let $\Sigma$ be a complete fan in $\R^2$.  By Proposition \ref{codim 1}, a maximal degree $d$ branched cover $\bf: \DDelta \rightarrow \SSigma$ is unramified along every cone $\sigma \neq 0_\Delta$ in $\Delta$.  It follows that $\bf|_{|\Delta|\smallsetminus \{ 0_\Delta \} }$ is a $d$-sheeted cover of $\R^2 \smallsetminus \{ 0 \}$, giving a partition $\lambda = (\lambda_1, \ldots, \lambda_s)$ of $d$, where $|\Delta| \smallsetminus \{ 0_\Delta \}$ has connected components $V_1, \ldots, V_s$ and the degree of the cover $V_i \rightarrow (\R^2 \smallsetminus \{ 0 \} )$ is $\lambda_i$.

Conversely, given a partition $\lambda = (\lambda_1, \ldots, \lambda_s)$ of $d$, there is a maximal branched cover $\SSigma_\lambda \rightarrow \Sigma$ constructed as follows.  For a positive integer $k$, let $V_k$ be the unique connected $k$-sheeted cover of $\R^2 \smallsetminus \{ 0 \}$, with
\[
|\Sigma_k| = V_k \sqcup \{ 0 \},
\]
with the obvious topology and projection to $\R^2$.  As topological spaces, $|\Sigma_k| \cong |\Sigma| \cong \C$, and the homeomorphisms may be chosen such that the projection is given by $z \mapsto z^k$.  Then $|\Sigma_k|$ inherits a natural cone complex structure ``pulled back" from $\SSigma$, where a cone in $\Sigma_k$ is the closure of a connected component of the preimage of the relative interior of a cone in $\Sigma$.  The cone complex $\SSigma_k$ is weighted by $\bw$, where $\bw(V_k) = 1$ and $\bw(0_{\Sigma_k}) = k$, and the projection to $\SSigma$ is a branched cover of degree $k$.  Then
\[
\SSigma_\lambda = \SSigma_{\lambda_1} \! \vee \cdots \vee \SSigma_{\lambda_s}
\]
is a maximal degree $d$ branched cover of $\SSigma$.
\end{example}

\begin{prop}
Let $\Sigma$ be a complete fan in $\R^3$, with $\bf: \DDelta \rightarrow \SSigma$ a branched cover.  Let $S \subset \R^3$ be the unit sphere, and let $B$ be the finite set of intersection points of $S$ with the rays of $\Sigma$.  Then $\bf$ is maximal if and only if $\bf^{-1}(S)$ is a topological manifold and $\bf^{-1}(x)$ consists of $d$ distinct points for every $x$ in $S \smallsetminus B$. 
\end{prop}

\begin{proof}
Suppose $\bf$ is maximal.  By Proposition \ref{codim 1}, $\bf$ is unramified away from the rays of $\Sigma$.  Therefore, $\bf|_{\bf^{-1}(S \smallsetminus B)}$ is a $d$-sheeted covering of $S \smallsetminus B$.  In particular, $\bf^{-1}(x)$ is $d$ distinct points for every $x$ in $S \smallsetminus B$, and to show that $\bf^{-1}(S)$ is a manifold it will suffice to check near a point $y \in \bf^{-1}(B)$.

Let $y$ be a point in $\bf^{-1}(B)$, let $U$ be a small open disc around $\bf(y)$ in $S$, and let $V$ be the connected component of $\bf^{-1}(U)$ containing $y$.  We must show that $V$ is homeomorphic to a disc.  We know that $V \smallsetminus \{y\}$ is a covering space of the annulus $U \smallsetminus \{ \bf(y) \}$, so $V$ is homeomorphic to a disc if and only if $V \smallsetminus \{ y \}$ is connected.  If $V \smallsetminus \{ y \}$ is disconnected, then we can pull apart the ramification along the ray in $\Delta$ containing $y$, replacing that ray with one copy for each connected component $V_i$ of $V \smallsetminus \{ y \}$ with weight equal to the degree of the covering $V_i \rightarrow U \smallsetminus \{ \bf(y) \}$.  Since $\bf$ is maximal, we conclude that $V \smallsetminus \{ y \}$ is connected, and hence $\bf^{-1}(S)$ is a manifold.

Conversely, suppose $\bf^{-1}(S)$ is a manifold and $\bf^{-1}(x)$ is $d$ distinct points for every $x$ in $S \smallsetminus B$.  Then $\bf$ is unramified away from the rays of $\Delta$.  If $\bf$ is not maximal, then we can pull apart some of the ramification along one of the rays in $\Delta$.  Since $\bf^{-1}(S)$ is a manifold, by hypothesis, for any $x \in B$ with $U$ a small open disc around $x$, $\bf^{-1}(x)$ contains exactly one point in the closure of each connected component of $\bf^{-1}(U \smallsetminus \{ x \} )$.  It follows that none of the ramification along the rays whose images contain $x$ can be pulled apart.  Therefore $\bf$ is maximal, as required.
\end{proof}

For a positive integer $d$, let $S_d$ be the symmetric group on $d$ letters.  

\begin{thm} \label{bijections}
Let $\Sigma$ be a complete fan in $\R^3$.  Let $S \subset \R^3$ be the unit sphere, and let $B$ be the finite set of intersection points of $S$ with the rays of $\Sigma$.  Then there are natural bijections between the following sets.
\begin{enumerate}

\item Isomorphism classes of maximal degree $d$ branched covers of $\Sigma$.

\item Isomorphism classes of degree $d$ holomorphic branched covers of the Riemann surface $\P^1(\C) \cong S$ whose branch locus is contained in $B$.

\item Isomorphism classes of $d$-sheeted covering spaces of $S \smallsetminus B$.

\item Conjugacy classes of representations $\pi_1(S \smallsetminus B) \rightarrow S_d$.
\end{enumerate}
\end{thm}

\begin{proof}
The map from (2) to (3) takes a branched cover $\pi : Y \rightarrow S$ to the covering space $\pi^{-1}(S \smallsetminus B)$, and the map from (3) to (4) takes a covering space to its monodromy representation.  That these are bijections is a standard result in the classical theory of Riemann surfaces.  See \cite[III.4]{Miranda95}.

To prove the theorem, we give a bijection between (1) and (3).  First, we give a map from branched covers of $\Sigma$ to covering spaces of $S \smallsetminus B$.  By Proposition \ref{top to poset}, part (2), a maximal branched cover $f: \Delta \rightarrow \Sigma$ induces a branched cover of cone complexes $\bf: \DDelta \rightarrow \SSigma$, which is maximal by Proposition \ref{top poset equivalence}.  By Proposition \ref{codim 1}, $\bf$ is unramified along codimension one cones.  Therefore $\bf|_{\bf^{-1} (S \smallsetminus B)}$ is a $d$-sheeted cover of $S \smallsetminus B$, and the map from (1) to (3) is given by $f \mapsto \bf|_{\bf^{-1}(S \smallsetminus B)}$.

Next, we give a map from covering spaces of $S \smallsetminus B$ to branched covers of $\Sigma$.  By the classical bijection from (3) to (2), a $d$-sheeted cover $\pi_0 : Y_0 \rightarrow S \smallsetminus B$ extends uniquely to a degree $d$ holomorphic branched cover $\pi: Y \rightarrow S$.  Then the cone over $Y$, denoted by
\[
CY \ = \ \bigl( Y \times [0,\infty) \bigr) \ / \ \bigl( Y \times \{ 0 \} \bigr)
\]
has a natural projection to $|\Sigma|$, the cone over $S$.  We give $CY$ the cone complex structure induced by this projection, where a cone in $CY$ is the closure of a connected component of the preimage of the relative interior of a cone in $\Sigma$.  Then $CY$ is weighted by $\bw$, where $\bw(y,t)$ is the ramification index of $\pi$ at $y$ for $t \neq 0$, and $\bw(0_{CY} ) = d$.  The projection from $C Y$ to $\SSigma$ is a branched cover of degree $d$, inducing a degree $d$ branched cover $f$ of $\Sigma$, which is maximal by Proposition \ref{top poset equivalence}.  The map from (3) to (1) is given by $\pi_0 \mapsto f$.

It is straightforward to check that the maps $f \mapsto \bf|_{\bf^{-1}(S \smallsetminus B)}$ and $\pi_0 \mapsto f$ are inverse to each other, up to natural isomorphisms, inducing a bijection on isomorphism classes, as required.
\end{proof}

\begin{remk}
The fundamental group of the complement of $m$ points in a 2-sphere is free on $(m-1)$ generators.  Therefore, if $\Sigma$ is a complete fan in $\R^3$ with $m$ rays, the bijection between (1) and (4) in Theorem \ref{bijections} implies that the number of isomorphism classes of maximal degree $d$ branched cover of $\Sigma$ is the number of $S_d$-conjugacy classes in $S_d^{(m-1)}$.  In particular, this number is bounded above by $d!^{(m-1)}$ and bounded below by $d!^{(m-2)}$.
\end{remk}

We close this section with a combinatorial problem with interesting implications in algebraic geometry. If $\Sigma$ is a complete fan with a degree $d$ branched cover $\bf: \DDelta \rightarrow \SSigma$, then a piecewise-linear function $\Psi \in \PL(\Delta)$ determines a multiset $\bu(\sigma) \subset M$, for each maximal cone $\sigma \in \Sigma$, where the multiplicity of $u$ in $\bu(\sigma)$ is the sum of the multiplicities of the cones $\sigma'$ in $\DDelta$ such that $f(\sigma') = \sigma$ and $\Psi|_{\sigma'} = \bf^* u$.  We say that $\Psi$ is trivial if there is some multiset $\bu \subset M$ such that $\bu(\sigma) = \bu$ for every maximal cone $\sigma \in \Sigma$.

\begin{prob} \label{problem}
Does every complete fan have a branched cover with a nontrivial piecewise-linear function on it?
\end{prob}

\noindent An affirmative answer to Totaro's Question would imply an affirmative solution to Problem~\ref{problem}.  In particular, any example of a complete fan with no nontrivial piecewise linear functions on any of its branched covers would give a counterexample to the resolution property.  On the other hand, an algorithm for producing a branched cover of a complete fan with a nontrivial piecewise linear function  might lead to a proof that every complete toric variety has a nontrivial vector bundle on it.  Either way, a solution to Problem~\ref{problem}, which is purely combinatorial, should be interesting to algebraic geometers.  

Section \ref{examples} gives examples of complete three-dimensional fans with no nontrivial piecewise linear functions on branched covers of degree two or three.  It is not known whether these fans have higher degree branched covers with nontrivial piecewise-linear functions on them.

\section{From toric vector bundles to branched covers of fans}

Let $X = X(\Delta)$ be a toric variety, and let $\E$ be a rank $r$ toric vector bundle on $X$.   For $\sigma \in \Delta$, recall that $M_\sigma = M / (\sigma^\perp \cap M)$ is the group of integral linear functions on $\sigma$.  

\subsection{The branched cover associated to a toric vector bundle} \label{branched covers and chern classes}

Any toric vector bundle on an affine toric variety $U_\sigma$ splits equivariantly as a sum of toric line bundles whose underlying line bundles are trivial.  For each $\sigma \in \Delta$, there is a unique multiset $\bu(\sigma) \subset M_\sigma$ such that 
\[
\E|_{U_\sigma} \cong  \bigoplus_{[u] \in \bu(\sigma)} \L_{[u]},
\]
where $\L_{[u]} = \O(\divisor \chi^u)$ for any lift $u \in M$ of $[u]$.  See \cite[Section 2]{moduli} for these and other basic properties of toric vector bundles.

Let $\Delta_\E$ be the poset
\[
\Delta_\E = \{ (\sigma, [u]) \ | \ \sigma \in \Delta \mbox{ and } [u] \in \bu(\sigma) \},
\]
where $(\sigma, [u] ) \geq (\tau, [u'])$ if and only if $\sigma \succeq \tau$ and $[u]|_\tau = [u']$.  Let $w: \Delta_\E \rightarrow \Z^+$ be the weight function such that $w((\sigma, [u]))$ is the multiplicity of $[u]$ in $\bu(\sigma)$.

\begin{prop}
The projection $(\Delta_\E, w) \rightarrow \Delta$ given by $(\sigma, [u]) \mapsto \sigma$ is a poset branched cover of degree $r$.
\end{prop}

\begin{proof}
For $\tau \preceq \sigma$, $\bu(\tau)$ is the multiset of restrictions of elements of $\bu(\sigma)$ to $\tau$, so the proposition follows directly from Definition~\ref{branched posets}.
\end{proof}

\begin{defn}
The weighted cone complex $(\DDelta_\E, \bw)$ is the fibered product
\[
\DDelta_\E = \DDelta \times_{\Delta} \Delta_\E,
\]
with weight function $\bw = p_2^* w$.  The branched cover associated to $\E$ is the projection
\[
p_1 : \DDelta_\E \rightarrow \DDelta.
\]
\end{defn}

\noindent The fact that $p_1 : (\DDelta_\E, \bw) \rightarrow \DDelta$ is a branched cover of degree $r$ follows directly from Proposition \ref{top to poset}.

\begin{defn}
The piecewise-linear function $\Psi_\E \in \PL(\DDelta_\E)$ associated to $\E$ is given by
\[
\Psi_\E|_{(\sigma, [u])} = p_1^*[u].
\]
\end{defn}

\noindent In particular, if $\E = \O(D)$ is a toric line bundle, then $\DDelta_\E = \DDelta$ and $\Psi_\E = \Psi_D$ is the usual piecewise-linear function associated to the $T$-Cartier divisor $D$.

\subsection{Equivariant Chern classes in terms of branched covers.}

The combinatorial data $(\Delta_\E, \Psi_\E)$ determine the equivariant total Chern class of $\E$ as follows.  The equivariant operational Chow cohomology ring $A^*_T(X)$ is canonically isomorphic to the ring of integral piecewise-polynomial functions 
\[
\PP^*(\DDelta) = \{ \Psi : |\Delta| \rightarrow \R \ | \ \Psi|_\sigma \in \Sym^* M_\sigma \mbox{ for all } \sigma \in \Delta \},
\]
and the $i$-th equivariant Chern class $c_i^T(\E)$ corresponds to the piecewise-polynomial function $\Psi^i_\E$, where
\[
\Psi^i_\E|_\sigma = e_i (\bu(\sigma)),
\]
the $i$-th elementary symmetric function in the multiset of linear functions $\bu(\sigma)$ \cite{chowcohom}.

We write $c_T(\E)$ for the equivariant total Chern class 
\[
c_T(\E) = 1 + c_1^T(\E) + \cdots + c_r^T(\E).
\]

\begin{prop} \label{chern classes}
Let $\E$ and $\E'$ be equivariant rank $r$ toric vector bundles on $X$.  Then $c_T(\E) = c_T(\E')$ if and only if there is an isomorphism of branched covers $g: \DDelta_{\E'} \rightarrow \DDelta_\E$ such that $g^*\Psi_\E = \Psi_{\E'}$.
\end{prop}

\begin{proof}
Say $\E'|_{U_\sigma} \cong  \bigoplus_{[u'] \in \bu'(\sigma)} \L_{[u']}$, for $\sigma \in \Delta$.  Then, by the basic properties of elementary symmetric functions, $e_i(\bu(\sigma)) = e_i(\bu'(\sigma))$ for all $1 \leq i \leq r$ if and only if $\bu(\sigma)  = \bu'(\sigma)$.  Therefore, $c_T(\E) = c_T(\E')$ if and only if $\bu(\sigma) = \bu'(\sigma)$ for every $\sigma \in \Delta$, and the proposition follows.
\end{proof}

\section{Complete toric varieties with few toric vector bundles}

Fix a complete toric variety $X = X(\Delta)$.  The category of toric vector bundles on $X$ is given by Klyachko's Classification Theorem \cite{Klyachko90}.

\begin{KCT}
The category of toric vector bundles on $X(\Delta)$ is naturally equivalent to the category of finite-dimensional $k$-vector spaces $E$ with collections of decreasing filtrations $\{E^\rho(i)\}$ indexed by the rays of $\Delta$, satisfying the following compatibility condition.

For each cone $\sigma \in \Delta$, there is a decomposition $E = \displaystyle{\bigoplus_{[u] \in M_\sigma}} E_{[u]}$ such that
\[
E^\rho(i) = \sum_{[u](v_\rho) \, \geq \, i} E_{[u]},
\]
 for every $\rho \preceq \sigma$ and $i \in \Z$.

\end{KCT}

\noindent The dimension of $E_{[u]}$ is the multiplicity of $[u]$ in the multiset of linear functions $\bu(\sigma) \subset M_\sigma$ associated to the corresponding toric vector bundle.

Natural operations on vector bundles, such as direct sums and tensor products, may be described easily in terms of Klyachko's filtrations.  In the proof of Theorem~\ref{trivial rank 3} below, the following characterization of the dual of a toric vector bundle will be helpful.

\begin{lem} \label{dual filtrations}
Let $\E$ be a toric vector bundle corresponding to the vector space with filtrations $(E, \{E^\rho(i)\})$.  The dual bundle $\E^\vee = \mathcal{H}om(\E, \O_X)$ corresponds to the dual vector space $E^\vee = Hom(E,k)$ with filtrations $(E^\vee){^\rho}(i) = E^\rho(1-i)^\perp$.
\end{lem}

\begin{proof}
It suffices to consider the case where $X = U_\rho$.  Then $\E$ splits as a sum of line bundles, so we may assume that $\E = \O(d D_\rho)$ for some integer $d.$  Then $\E^\vee = \O(-dD_\rho)$, so
\[
E^\rho(i) = \left\{ \begin{array}{ll} E & \mbox{ for } i \leq d, \\ 
						 0 & \mbox{ for } i > d, \end{array} \right.
						 \mbox{ \ \ \ \ and \ \ \ \ }
(E^\vee)^\rho(i) = \left\{ \begin{array}{ll} E^\vee & \mbox{ for } i \leq -d, \\
   							    0 & \mbox{ for } i > -d, \end{array} \right.
\]							    					
and the lemma follows.
\end{proof}

\subsection{Toric vector bundles with trivial Chern class} 

We say that a toric vector bundle $\E$ on $X$ has trivial Chern class if the equivariant total Chern class $c_T(\E)$ in $A^*_T(X) \cong \PP^*(\Delta)$ is given by a globally polynomial function.  By the basic properties of elementary symmetric functions, this is equivalent to saying that $\E$ has trivial Chern class if there is a multiset $\bu \subset M$ such that $\bu(\sigma) = \bu$ for every maximal cone $\sigma$ in $\Delta$.  The toric vector bundles whose underlying bundles are trivial are exactly those isomorphic to some $\bigoplus_{u \in \bu} \O_X(\divisor \chi^{u})$, so what we are calling trivial Chern classes are exactly the equivariant Chern classes of the toric vector bundles on $X$ whose underlying bundles are trivial. The following are examples of toric vector bundles with trivial Chern class whose underlying vector bundles are not trivial.

\begin{example}
Suppose $N = \Z$ and $\Delta$ is the unique complete fan in $N_\R = \R$.  Then $X(\Delta)$ is isomorphic to $\P^1$.  Let $\rho$ and $\rho'$ be the rays of $\Delta$, spanned by $e_1$ and $-e_1$, respectively.  Let $E$ be a two-dimensional $k$-vector space, with $L \subset E$ a one-dimensional subspace.  Let $\E$ be the toric vector bundle on $X$ given by the vector space $E$ with filtrations
\[
E^\rho(i) = \left\{ \begin{array}{ll} E & \mbox{for $i < 1$,} \\
						  L & \mbox{for $i = 1$,} \\
						 0 & \mbox{for $i > 1$,}
						  \end{array} \right.
\mbox{ \ \ \ \ and \ \ \ \ \ \ }
E^{\rho'}(i) = \left\{ \begin{array}{ll}  E & \mbox{for $i < 0$}, \\
						  L & \mbox{for $i = 0$,} \\
						0 & \mbox{for $i > 0$.}
						  \end{array} \right.
\]
Then $\bu(\rho) = \bu(\rho') = \{ 0, e_1^* \}$, so $\E$ has trivial Chern class.  However, $\E$ is not trivial.  The underlying bundle of $\E$ is isomorphic to $\O(1) \oplus \O(-1)$, where the subbundle $\O(1) \subset \E$ has fiber $L$ over $x_0$.
\end{example}

\begin{example} \label{rank 3 example}
Suppose $N = \Z^2$ and $\Delta$ is the complete fan in $N_\R = \R^2$ whose rays $\rho_1$, $\rho_2$, and $\rho_3$ are spanned by $e_1$, $e_2$, and $(-e_1-e_2)$, respectively.  Then $X(\Delta)$ is isomorphic to $\P^2$.  Let $E$ be a three-dimensional $k$-vector space, with $F \subset E$ a two dimensional subspace.  Let $L_1$, $L_2,$ and $L_3$ be distinct one-dimensional subspaces of $F$, and let $\E$ be the toric vector bundle on $X$ given by the vector space $E$ with filtrations
\[
E^{\rho_1}(i) = \left\{ \begin{array}{ll} E & \mbox{for $i < 1$,} \\
						  L_1 & \mbox{for $i = 1$,} \\
						 0 & \mbox{for $i > 1$,}
						  \end{array} \right.
\mbox{ \ \ \ \ \ \ \ }
E^{\rho_2}(i) = \left\{ \begin{array}{ll} E & \mbox{for $i < 1$,} \\
						  L_2 & \mbox{for $i = 1$,} \\
						  0 & \mbox{for $i > 1$,}
						  \end{array} \right.
\]
and
\[
E^{\rho_3}(i) =  \left\{ \begin{array}{ll} E & \mbox{for $i < 0$,}\\
						  L_3 & \mbox{for $i = 0$,} \\
						 0 & \mbox{for $i > 0$.} 
						  \end{array} \right.
\]
Then $\bu(\sigma) = \{0, e_1^*, e_2^*\}$ for each maximal cone $\sigma$ in $\Delta$, so $\E$ has trivial Chern class.  However, $\E$ is not trivial.  There is a split surjection $\E \rightarrow \O(-1)$ with kernel $\F \cong T\P^2 \otimes \O(-1)$, where $\F \subset \E$ is the toric subbundle with fiber $F$ over $x_0$.
\end{example}

\noindent One can construct similar examples of nontrivial rank $(n+1)$ toric vector bundles with trivial Chern class on $\P^n$, for any $n$, and, more generally, on any complete $n$-dimensional $\Q$-factorial toric variety with Picard number 1.				

\begin{thm} \label{trivial rank 2}
Let $X$ be a complete toric variety and let $\E$ be a nontrivial toric vector bundle of rank two on $X$ with trivial Chern class.   Then there is a surjective toric morphism $\pi: X \rightarrow \P^1$ and a nontrivial toric vector bundle $\F$ on $\P^1$ such that $\E \otimes \O(\divisor \chi^u)$ is equivariantly isomorphic to $\pi^* \F$, for some $u \in M$.
\end{thm}

\noindent In the situation of Proposition \ref{trivial rank 2}, the underlying bundle of $\F$ is necessarily isomorphic to $\O(a) \oplus \O(-a)$, for some integer $a$.  

The classification of rank three toric vector bundles with trivial Chern class is similar, but slightly more complicated; such bundles are pulled back from either $\P^1$ or a surface with small Picard number.

\begin{thm} \label{trivial rank 3}
Let $X$ be a complete toric variety and let $\E$ be a nontrivial toric vector bundle of rank three on $X$ with trivial Chern class.   Then there is a surjective toric morphism $\pi: X \rightarrow Y$ and a nontrivial toric vector bundle $\F$ on $Y$ such that $\E \otimes \O(\divisor \chi^u)$ is equivariantly isomorphic to $\pi^* \F$, for some $u \in M$, where $Y$ is isomorphic to either $\P^1$ or a projective toric surface with Picard number one, two, or four.
\end{thm}

\begin{remk}
In particular, any complete toric variety that has a nontrivial toric vector bundle with trivial Chern class of rank less than or equal to three admits a nonconstant morphism to a projective variety. Since a toric variety with no nontrivial nef line bundles admits no such morphisms, any toric vector bundle of rank less than or equal to three with trivial Chern class on such a variety is trivial.  Varieties with no nontrivial nef line bundles may be singular and have no nontrivial line bundles at all, like the examples in Section \ref{examples}, or they may be smooth, like the examples studied in \cite{nonef}.
\end{remk}

In the proofs of Theorems \ref{trivial rank 2} and \ref{trivial rank 3}, we use collections of decreasing $\R$-graded filtrations $\{ E^v(t) \}$ indexed by points $v \in |\Delta|$ that, roughly speaking, continuously interpolate the filtrations $E^\rho(i)$ from Klyachko's classification.

\subsection{Continuously interpolating Klyachko's filtrations}

Let $v$ be a point in the relative interior of a cone $\sigma$ in $\Delta$.  Recall that $E^\sigma_u$ is the image in the fiber $E = \E_{x_0}$ of the space of $\chi^u$-isotypical sections $\Gamma(U_\sigma, \E)_u$.  For a real number $t$, we define
\[
E^v(t) \ = \sum_{\<u, v \> \, \geq \, t } E^\sigma_u.
\]
The filtrations $\{E^v(t) \}$ are a straightforward generalization of the filtrations $E^\rho(i)$; if $i$ is an integer and $v_\rho$ the primitive generator of a ray $\rho$ in $\Delta$, then
\[
E^{v_\rho}(i) = E^\rho(i).
\]

Let $J(v) \subset \{1, \ldots, r \}$ be the set of dimensions of the nonzero vector spaces $E^v(t)$, for $t \in \R$, and let $\Fl(v) \in \FFl_{J(v)}(E)$ be the partial flag consisting of exactly these subspaces.  Fix a cone $\sigma \in \Delta$, with a multiset $\bu(\sigma) \subset M_\sigma$ and a splitting $E = \bigoplus_{[u] \in \bu(\sigma)} E_{[u]}$ verifying the compatibility condition in Klyachko's Classification Theorem.  Then it follows directly from the definitions that, for $v \in \sigma$,
\[
E^v(t) \ = \sum_{[u](v) \, \geq \, t} E_{[u]}.
\]
In particular, the flag $\Fl(v)$ depends only on the preorder on $\bu(\sigma)$ where $[u] \leq [u']$ if and only if $[u](v) \leq [u'](v)$.  Since there are only finitely many preorders on $\bu(\sigma)$, there are only finitely many partial flags that appear as $\Fl(v)$, for $v$ in $\sigma$.  Let $\Po(\sigma)$ be this finite set of partial flags, with the partial ordering where $x \leq y$ in $\Po(\sigma)$ if and only if every subspace appearing in $x$ also appears in $y$.  Give $\Po(\sigma)$ the poset topology, as in Section \ref{poset section}, where
\[
\Po(\sigma)_x = \{ y \in \Po(\sigma) \ | \ y \leq x \} \mbox{ \ and \ } \Po(\sigma)^x = \{ y \in \Po(\sigma) \  | \ y \geq x \}
\]
are basic closed sets and open sets, respectively, for $x \in \Po(\sigma)$.
						  
\begin{lem}
The map $f_\sigma: \sigma \rightarrow \Po(\sigma)$ given by $v \mapsto \Fl(v)$ is continuous.
\end{lem}

\begin{proof}
A partial flag $x$ in $\Po(\sigma)$ determines a unique preorder $\leq_x$ on $\bu(\sigma)$ such that $\Fl(v) \leq x$ if and only if $[u](v) \leq [u'](v)$ whenever $[u] \leq_x [u']$.  Therefore, the preimage under $f_\sigma$ of the basic closed set $\Po(\sigma)_x$ is 
\[
f_\sigma^{-1}(\Po(\sigma)_x) = \{ v \in \sigma \ | \ [u](v) \leq [u'](v) \mbox{ for } [u] \leq_x [u'] \},
\]
which is a closed polyhedral cone.  In particular, the preimage of a basic closed set is closed, so $f_\sigma$ is continuous.
\end{proof}

The maps $f_\sigma$, for $\sigma \in \Delta$, glue together to give a continuous map $f_\Delta : |\Delta| \rightarrow \Po(\Delta)$, where $\Po(\Delta)$ is the partially ordered set of partial flags that appear as $\Fl(v)$ for $v \in |\Delta|$, and $\Po(\Delta)$ has the poset topology.  For $j \in J(v)$, we write $\Fl(v)_j$ for the $j$-dimensional subspace appearing in $\Fl(v)$.

\begin{prop} \label{constant on components}
The subspace $\Fl(v)_j$ of $E$ is constant on each connected component of the set of points $v$ in $|\Delta|$ such that $J(v)$ contains $j$.
\end{prop}

\begin{proof}
Let $F$ be a $j$-dimensional subspace of $E$.  By the continuity of $f_\Delta$, the set of points $v$ in $|\Delta|$ such that $\Fl(v)$ contains $F$ is both open and closed in the locus where $J(v)$ contains $j$.
\end{proof}

\begin{proof}[Proof of Theorem \ref{trivial rank 2}]
Say $\bu \subset M$ is the multiset such that $\E|_{U_\sigma} \cong \bigoplus_{u \in \bu} \L_u$ for every maximal cone $\sigma$ in $\Delta$.  If $\bu$ consists of one element with multiplicity two, then $\Fl(v) = E$ for all $v$ in $N_\R$, and it is easy to see that $\E$ is trivial.

So assume $\bu$ has two distinct elements, $u_1$ and $u_2$.  By Proposition \ref{constant on components}, $\Fl(v)$ is constant on the open halfspaces $U_1$ and $U_2$, where $u_1 > u_2$ and where $u_2 > u_1$, respectively.  Let $L_1$ and $L_2$ be the one-dimensional subspaces of $E$ such that $\Fl(v)_1 = L_i$ for $v$ in $U_i$.  Then $\E$ is trivial if and only if $L_1 + L_2 = E$, so it suffices to consider the case when $L_1 = L_2$.

Suppose $L_1 = L_2 = L$, and choose $L' \subset E$ such that $E$ splits as $E \cong L \oplus L'$.  It follows from Klyachko's Classification Theorem that $\E$ splits as a sum of line bundles $\E \cong \L \oplus \L'$, where $L$ and $L'$ are the fibers over $x_0$ of $\L$ and $\L'$, respectively.  The corresponding piecewise-linear functions are given by
\[
\Psi_\L (v) = \max\{ \< u_1, v\> , \<u_2, v \> \} \mbox{ \ \ and \ \ } \Psi_{\L'}(v) = \min\{ \< u_1, v\> , \<u_2, v \> \}.
\]
Therefore, every maximal cone of $\Delta$ is contained in either $\overline{U}_1$ or $\overline{U}_2$, the maximal domains of linearity of $\Psi_\L$ and $\Psi_\L'$.  Let $N_0$ be the kernel of $(u_1 - u_2)$, so there is a projection $X \rightarrow \P^1$ corresponding to the quotient of lattices $N \rightarrow N/ N_0$.  Since $\Psi_\L - u_1$ and $\Psi_{\L'} - u_1$ are the pullbacks of integral piecewise-linear functions on $(N/N_0)_\R$, the twist $\E \otimes \O(\divisor \chi^{-u_1})$ is equivariantly isomorphic to the pullback of a toric vector bundle on $\P^1$.
\end{proof}

The proof of Theorem \ref{trivial rank 3} is somewhat more complicated than the proof of Theorem \ref{trivial rank 2}, in part because rank three toric vector bundles with trivial Chern class do not necessarily split as a sum of line bundles.  For instance, the bundle in Example \ref{rank 3 example} has an irreducible rank two subbundle with fiber $F$ over $x_0$.  The following proposition describes the pullbacks of toric vector bundles under toric morphisms in terms of the filtrations $\{ E^v (t) \}$.  Let $\Delta'$ be a fan in $N'_\R$, and let $\varphi: N' \rightarrow N$ be a map of lattices inducing a toric morphism $f: X(\Delta') \rightarrow X$.  Then $f$ maps the identity $x_0'$ in the dense torus $T' \subset X(\Delta')$ to $x_0$, so the fiber of $f^* \E$ over $x_0'$ is canonically identified with $E$.

\begin{prop} \label{pullback filtrations}
The filtrations determining $f^*\E$ are given by $E^{v'}(t) = E^{\varphi(v')} (t)$, for $v' \in N'_\R$.
\end{prop}

\begin{proof}
Let $M' = \Hom(N', \Z)$ be the character lattice of $T'$, and  let $\sigma'$ be a cone in $\Delta'$ such that $\varphi(\sigma')$ is contained in $\sigma \in \Delta$.  Then $f$ maps $U_\sigma'$ into $U_\sigma$, so there is a $M'$-graded isomorphism $\Gamma(U_\sigma', f^*\E) \cong \Gamma(U_\sigma, \E) \otimes_{k[U_\sigma]} k[U_{\sigma'}],$ and it follows that
\[
E^{\sigma'}_{u'} = \sum_{(\varphi^*u - u') \in (\sigma')^\vee } E^\sigma_u,
\]
for $u' \in M'$.  Therefore,
\[
E^{v'}(t) \ = \sum_{\<u', v \> \, \geq \, t } E^{\sigma'}_{u'} \ = \sum_{\< \varphi^*u, v' \> \, \geq \, t } E^\sigma_u = \sum_{\< u, \varphi(v') \> \, \geq \, t} E^\sigma_u,
\]
which is equal to $E^{\varphi(v')}(t)$ by definition.
\end{proof}

\begin{proof}[Proof of Theorem \ref{trivial rank 3}]
Let $\bu \subset M$ be the multiset such that $\E|_{U_\sigma} \cong \bigoplus_{u \in \bu} \L_u$ for every maximal cone $\sigma$ in $\Delta$.  Say $\bu = \{u_1, u_2, u_3 \}$.  We first handle the cases where the $u_i$ are collinear in $M_\R$.

If the $u_i$ are all the same, then $\Fl(v) = E$ for all $v \in N_R$, and $\E$ is trivial.  If $u_1$ and $u_2$ are distinct, and $u_3$ lies on the line spanned by $u_1$ and $u_2$, then $\Fl(v)$ is constant on the halfspaces where $u_1 > u_2$ and where $u_2 > u_1$.  Arguments similar to those in the proof of Theorem \ref{trivial rank 2} then show that either $\E$ is trivial, or the projection $N \rightarrow N/N_0$ induces a toric morphism $X \rightarrow \P^1$ such that $\E \otimes \O(\divisor \chi^{u_1})$ is equivariantly isomorphic to the pullback of a toric vector bundle on $\P^1$.

It remains to consider the cases where $u_1, u_2$, and $u_3$ are not collinear.  Let
\[
U_i = \{ v \in N_\R \ | \ \< u_i, v \> < \<u_j, v \> \mbox{ for } j \neq i \},
\]
and
\[
U^k = \{v \in N_\R \ | \ \<u_j, v \> < \< u_k, v \> \mbox{ for } j \neq k \}.
\]
Let $U_{ijk}$ be the intersection
\[
U_{ijk} = U_i \cap U^k = \{ v \in N_\R \ | \ \<u_i, v \> < \< u_j, v \> < \<u_k, v \> \}.
\]
The configuration of these regions is illustrated below.

\vspace{10 pt}

\begin{center}
\begin{picture}(340,110)(0,-10)

\put(0,0){\line(1,1){40}}
\put(40,40){\line(0,1){50}}
\put(40,40){\line(1,0){50}}

\put(65,65){\makebox(0,0)[l]{$U^3$}}
\put(50,10){\makebox(0,0)[l]{$U^1$}}
\put(10,50){\makebox(0,0)[r]{$U^2$}}

\put(130,0){\line(1,1){40}}
\put(170,40){\line(0,1){50}}
\put(170,40){\line(1,0){50}}
\put(170,40){\line(-1,0){50}}
\put(170,40){\line(0,-1){50}}
\put(170,40){\line(1,1){40}}

\put(195,10){\makebox(0,0)[l]{$U_{231}$}}
\put(155,-5){\makebox(0,0)[c]{$U_{321}$}}
\put(127,23){\makebox(0,0)[c]{$U_{312}$}}
\put(200,55){\makebox(0,0)[l]{$U_{213}$}}
\put(180,82){\makebox(0,0)[l]{$U_{123}$}}
\put(148,70){\makebox(0,0)[r]{$U_{132}$}}

\put(300,40){\line(-1,0){50}}
\put(300,40){\line(0,-1){50}}
\put(300,40){\line(1,1){40}}

\put(280,10){\makebox(0,0)[r]{$U_3$}}
\put(295,65){\makebox(0,0)[r]{$U_1$}}
\put(325,30){\makebox(0,0)[l]{$U_2$}}

\end{picture}
\end{center}

\vspace{10 pt}

\noindent By Proposition \ref{constant on components}, $\Fl(v)_2$ is constant on each $U_i$ and $\Fl(v)_1$ is constant on each $U^k$.  Let 
\[
L_k = \Fl(v)_1 \mbox{ for } v \in U^k,
\]
and
\[
V_i = \Fl(v)_2 \mbox{ for } v \in U_i.
\]
Then $L_k$ is contained in $V_i$ for $i \neq k$, and $\E$ is trivial if and only if $L_1 + L_2 + L_3 = E$.  So assume $L_1 + L_2 + L_3 \neq E$.  We complete the proof of the theorem case by case, according to which of the $L_k$ are equal, and which of the $V_i$ are equal.  Let $N_0 \subset N$ be the sublattice
\[
N_0 = \{ v \in N \ | \ \<u_1, v \> = \<u_2, v \> = \<u_3, v\> \}.
\]

\emph{Case 1.} Suppose $L_1 = L_2 = L_3 = L.$  Then there is a toric line bundle $\L \subset \E$ with fiber $L$ over $x_0$.  The piecewise-linear function $\Psi_\L$ is given by $\Psi_\L(v) = \max \{ \<u_1, v\> , \<u_2, v\>, \<u_3, v\> \}$, so every cone of $\Delta$ is contained in some $\overline{U}{^{k}}$.

\emph{Subcase 1.1.} Suppose $V_1 = V_2 = V_3 = V.$  Choose $L'$ and $L''$ such that $V \cong L \oplus L'$ and $E \cong L \oplus L' \oplus L''$.  Then $\E$ splits as a sum of line bundles $\E \cong \L \oplus \L' \oplus \L''$, where $L'$ and $L''$ are the fibers over $x_0$ of $\L'$ and $\L''$, respectively.
The corresponding piecewise-linear functions $\Psi_{\L'}$ are given by 
\[
\Psi_{\L'} (v) = \<u_j, v\> \mbox{ for } v \in \overline{U}_{ijk},
\]
and
\[
\Psi_{\L''}(v) = \min \{ \<u_1, v\> , \<u_2, v\>, \<u_3, v\> \}.
\]
Therefore, each cone of $\Delta$ must be contained in some $\overline{U}_{ijk}$. The projection $N \rightarrow N/N_0$ induces a surjective toric morphism $X \rightarrow X(\Sigma)$, where $\Sigma$ is the fan in $(N/N_0)_\R$ whose maximal cones are the images of the $\overline{U}_{ijk}$.  Now $X(\Sigma)$ is a projective toric surface with Picard number four, and $\E \otimes \O(\divisor \chi^{-u_1})$ is isomorphic to the pullback of a toric vector bundle on $X(\Sigma)$.

\emph{Subcase 1.2.} Suppose $V_1 = V_2  \neq V_3$.  Choose $L'$ and $L''$ such that $V_1 \cong L \oplus L'$ and $V_3 = L \oplus L''$.  Then $\E$ splits as a sum of line bundles $\E \cong \L \oplus \L' \oplus \L''$, where $L'$ and $L''$ are the fibers over $x_0$ of $\L'$ and $\L''$, respectively.  The corresponding piecewise-linear functions are given by
\[
\Psi_{\L'}(v)= \left\{ \begin{array}{l} \<u_j, v\> \mbox{ for } v \in \overline{U}_{ijk} \mbox{ and } j \neq 3, \\
                                                          \<u_i, v\> \mbox{ for } v \in \overline{U}_{ijk} \mbox{ and } j = 3, 
                                                          \end{array} \right.
\]                                         
and
\[
\Psi_{\L''}(v)= \left\{ \begin{array}{l} \<u_i, v\> \mbox{ for } v \in \overline{U}_{ijk} \mbox{ and } j \neq 3, \\
                                                          \<u_j, v\> \mbox{ for } v \in \overline{U}_{ijk} \mbox{ and } j = 3. 
                                                          \end{array} \right.
\]       
It follows that each cone of $\Delta$ is contained in either $\overline{U}{^1}$, $\overline{U}{^2}$, $\overline{U}_{123},$ or $\overline{U}_{213}$.  The projection $N \rightarrow N/N_0$ induces a surjective toric morphism $X \rightarrow X(\Sigma)$, where $\Sigma$ is the fan in $(N/N_0)_\R$ whose maximal cones are the images of $\overline{U}{^1}$, $\overline{U}{^2}$, $\overline{U}_{123},$ and $\overline{U}_{213}$.  Then $X(\Sigma)$ is a projective toric surface with Picard number two, and $\E \otimes \O(\divisor \chi^{-u_1})$ is pulled back from $X(\Sigma)$.

\emph{Subcase 1.3.} Suppose $V_1$, $V_2,$ and $V_3$ are distinct.  Then $\E / \L$ is irreducible, so $\E$ does not split as a sum of line bundles.  

Fix $k$, let $L'_k = L$, and choose one-dimensional subspaces $L'_i$ and $L'_j$ of $E$ such that $V_i$ and $V_j$ are equal to $L + L_i'$ and $L + L_j'$, respectively.  Then, for $v \in \overline{U}{^k}$,
\[
E^v(t) = \sum_{ \<u_i, v\> \geq t } L'_i.
\]
Let $\F = \E \otimes \O(\divisor \chi^{-u_1})$, so $\F$ corresponds to the vector space $F = E$ with filtrations
\[
F^v(t) = E^v(t + \<u_1, v\>),
\]
and let $\Sigma$ be the fan in $(N/N_0)_\R$ whose maximal cones are the images of the $\overline{U}{^k}$.  Then $X(\Sigma)$ is a projective toric surface with Picard number one.  Note that $F^v(t) = F^{v'}(t)$ whenever $v$ and $v'$ have the same image in $(N/N_0)_\R$, so the filtrations $F^v(t)$ descend to filtrations $F^{\overline{v}}(t)$, for $\overline{v}$ in $(N/N_0)_\R$.  These induced filtrations satisfy the compatibility condition in Klyachko's Classification Theorem for each cone in $\Sigma$, and hence define a toric vector bundle $\mathcal{G}$ on $X(\Sigma)$.  It then follows from Proposition~\ref{pullback filtrations} that $\F$ is isomorphic to $\pi^* \mathcal{G}$, where $\pi$ is the surjective toric morphism $X \rightarrow X(\Sigma)$ induced by the projection $N \rightarrow (N/N_0)_\R$.

\emph{Case 2.}  Suppose $L_1 = L_2 \neq L_3$.  Let $L = L_1$, $L' = L_3$, and $V = L + L'$.  Then $V_1$ and $V_2$ are both equal to $V$.

\emph{Subcase 2.1.} Suppose $V_3 = V$.  To prove the theorem for $\E$ it suffices to prove it for $\E^\vee$.  For $\E^\vee$, the theorem follows from Lemma \ref{dual filtrations} and Subcase 1.2.

\emph{Subcase 2.2.} Suppose $V_3 \neq V$.  Choose $L'' \subset E$ such that $E \cong L \oplus L' \oplus L''$.  Then $\E$ splits as a sum of line bundles $\E \cong \L \oplus \L' \oplus \L''$, where $L$, $L'$, and $L''$ are the fibers over $x_0$ of $\L$, $\L'$, and $\L''$, respectively.  The corresponding piecewise-linear functions are given by
\[
\Psi_\L(v) = \max \{ \<u_1, v\>, \<u_2, v\> \}, \mbox{\ \ \ \ \ }
\Psi_{\L''}(v) = \min \{ u_1, v\> , \<u_2, v \> \},
\]
and $\Psi_{\L'}(v) = \<u_3, v\>$.  It follows that the projection $N \rightarrow N/N'$ induces a morphism $X \rightarrow \P^1$, where $N'$ is the kernel of $u_1 - u_2$, such that $\E \otimes \O(\divisor \chi^{-u_1})$ is pulled back from $\P^1$.

\emph{Case 3.}  Suppose $L_1, L_2,$ and $L_3$ are distinct.  Then $V_1$, $V_2$, and $V_3$ are all equal to the two-dimensional subspace $L_1 + L_2 + L_3 \subset E$.  To prove the theorem for $\E$ it suffices to prove it for $\E^\vee$, and for $\E^\vee$, the theorem follows from Lemma \ref{dual filtrations} and Subcase 1.3.
\end{proof}

\subsection{Examples} \label{examples}

In this section, we study toric vector bundles of ranks two and three on specific complete toric threefolds that have no nontrivial line bundles.

\begin{example}[Eikelberg's threefold]  Let $N= \Z^3$ and let $\Sigma$ be the fan in $N_\R$ whose rays are generated by
\[
\begin{array}{lll}
v_1 = (-1,0,1), & v_2 = (0,-1,1), & v_3 = (1,2,1), \\
v_4 = (-1,0,-1), & v_3 = (0,-1,-1), & v_6 = (1,1,-1),
\end{array}
\]
and whose maximal cones are
\[
\begin{array}{lll}
\sigma_1 = \<v_1, v_2, v_3\>,  & \sigma_2 = \<v_4,v_5,v_6\>, & \sigma_3 = \<v_1,v_2,v_4,v_5\>, \\
 \sigma_4 = \<v_2,v_3,v_5,v_6\>, & \sigma_5 = \<v_1, v_3, v_4, v_6\>.
\end{array}
\]
Then $X(\Sigma)$ is complete, and Eikelberg showed that $\Sigma$ has no nontrivial piecewise-linear functions, and hence that $X =X(\Sigma)$ has no nontrivial line bundles \cite[Example~3.5]{Eikelberg92}.  

We claim that there do exist nontrivial rank two toric vector bundles on $X$.  By Theorem~\ref{trivial rank 2}, any such bundle must have a nontrivial equivariant Chern class, corresponding to a nontrivial piecewise-linear function on a branched cover of $\SSigma$.  Consider the unique maximal degree two branched cover $\DDelta \rightarrow \SSigma$ that is branched over exactly $\rho_1$ and $\rho_6$.  Let $\Psi$ be the piecewise-linear function on $\DDelta$ given by the pairs of linear functions
\[
\begin{array}{ll}
\bu(\sigma_1) = \{ (15,-15,3), (3,3,-9) \}, & \bu(\sigma_2) = \{ (16,-14,-4), (2,2,-2) \}, \\
\bu(\sigma_3) = \{ (12,-18,0), (6,6,-6)) \}, & \bu(\sigma_4) = \{ (24,-18,0), (-6,6,-6) \}, \\
 \bu(\sigma_5) = \{ (12,-6,0), (6,-6,-6) \}.
\end{array}
\]
Let $E$ be a two-dimensional $k$-vector space with distinct one-dimensional subspaces $L$, $L'$, and $L''$.  The following filtrations give a rank two toric vector bundle $\E$ on $X$ whose equivariant total Chern class $c_T(\E)$ is equal to $c(\Delta, \Psi)$.  The filtrations are
\[
\begin{array}{ll}
{ E^{\rho_1}(i) = \left\{ \begin{array}{ll} E & \mbox{ for } i \leq -12, \\ 0 & \mbox{ for } i > -12, \end{array} \right. }
& 
{ E^{\rho_2}(i) = \left\{ \begin{array}{ll} E & \mbox{ for } i  \leq -12,  \\ L & \mbox{ for } -12  <  i   \leq 18, \\ 0 & \mbox{ for } i  > 18, \end{array} \right. } \\ \\ 
 { E^{\rho_3}(i) = \left\{ \begin{array}{ll} E & \mbox{ for } i \leq -12, \\ L' & \mbox{ for } -12 < i \leq 0, \\ 0 &  \mbox{ for } i > 0, \end{array} \right. }
  &
  { E^{\rho_4}(i) = \left\{ \begin{array}{ll} E & \mbox{ for } i  \leq -12,  \\ L'' & \mbox{ for } -12  <  i  \leq 0 , \\ 0 & \mbox{ for } i  > 0, \end{array} \right. } \\ \\
 { E^{\rho_5}(i) = \left\{ \begin{array}{ll} E & \mbox{ for } i \leq 0, \\ L & \mbox{ for } 0 < i \leq 18, \\ 0 &  \mbox{ for } i > 18, \end{array} \right. }
  &
  { E^{\rho_6}(i) = \left\{ \begin{array}{ll} E & \mbox{ for } i \leq 6, \\ 0 & \mbox{ for } i > 6. \end{array} \right. }
  \end{array}
\]
It is straightforward to verify that these filtrations satisfy the compatibility condition in Klyachko's Classification Theorem, with $\bu(\sigma_i)$ as above, and the claim follows.
\end{example}

\begin{example}[Fulton's threefold] \label{Fulton's threefold}
Let $N = \Z^3$, and let $\Sigma$ be the fan in $N_\R$ whose rays are generated by 
\[
\begin{array}{llll}
v_1 = (1,2,3), & v_2 = (1,-1,1), & v_3 = (-1,-1,1), & v_4 = (-1,1,1), \\
v_5 = (1,1,-1), & v_6 = (1,-1,-1), & v_7 = (-1,-1,-1), & v_8 = (-1,1,-1),
\end{array}
\]
and whose maximal cones are
\[
\begin{array}{lll}
\sigma_1 = \< v_1,v_2, v_3, v_4\>, & \sigma_2 = \<v_5, v_6, v_7, v_8\>, & \sigma_3 = \< v_1, v_2, v_5, v_6 \>, \\
\sigma_4 = \< v_2, v_3, v_6, v_7\>, & \sigma_5 = \< v_3, v_4, v_7, v_8 \>, & \sigma_6 = \<v_1, v_4, v_5, v_8 \>.
\end{array}
\]
The toric variety $X = X(\Sigma)$ is complete, and Fulton has shown that $X$ has no nontrivial line bundles \cite[pp.\ 25--26, 72]{Fulton93}.

We now show that $X$ has no nontrivial rank two toric vector vector bundles.  Since $X$ has no nontrivial line bundles, there are no nonconstant morphisms from $X$ to projective varieties.  Therefore, it follows from Theorem \ref{trivial rank 2} that there are no nontrivial rank two toric vector bundles with trivial Chern class on $X$.  Hence, to show that there are no nontrivial rank two toric vector bundles on $X$, it is enough to show that there are no nontrivial piecewise-linear functions on maximal degree two branched covers of $\SSigma$.

Suppose $\bf:\DDelta \rightarrow \SSigma$ is a maximal degree two branched cover, with $\Psi \in \PL(\DDelta)$ a nontrivial piecewise-linear function.  Let $S \subset N_\R = \R^3$ be the unit sphere, and let $B$ be the finite set of intersection points of $S$ with the rays of $\Sigma$.  By Theorem \ref{bijections}, $\bf$ is uniquely determined by the monodromy representation
\[
\pi_1(S \smallsetminus B) \rightarrow S_2.
\]
The fundamental group $\pi_1(S \smallsetminus B)$ is generated by small loops around each of the points in $B$, so the monodromy representation is completely determined by the subset of these generators that map to the nontrivial element in $S_2$, which correspond exactly to the rays of $\Sigma$ over which $\bf$ is ramified.  The product of these generators is the identity, so the number of rays in this subset must be even.

We claim that it suffices to consider the case where no two rays over which $\bf$ is ramified span a two-dimensional cone.  Suppose, to the contrary, that $\bf$ is ramified over two rays $\rho$ and $\rho'$ that span a two-dimensional cone $\tau$ in $\Sigma$.  Let $\overline{\DDelta}$ be the branched cover of $\Sigma$ obtained by identifying the two preimages of $\tau$ in $\Delta$.  Since there is a single pair of rays in $\Delta$ that span both copies of $\tau$ in $\Delta$, there is a single linear function $[u] \in M_\tau$ such that $\Psi$ restricts to $\bf^*[u]$ on each of the two copies of $\tau$ in $\Delta$, and hence $\Psi$ is the pullback of a piecewise-linear function $\overline{\Psi} \in \PL(\overline{\Delta})$.  Let $\bf' : \DDelta' \rightarrow \SSigma$ be the maximal degree two branched cover that is ramified over all of the rays over which $\bf$ is ramified except $\rho$ and $\rho'$.  Then $\bf'$ factors through $\bg: \DDelta' \rightarrow \overline{\DDelta}$, and $\bg^* \overline{\Psi}$ is a nontrivial piecewise-linear function on $\DDelta'$.  The claim follows by induction on the number of rays over which $\bf$ is ramified.

There are eighteen maximal degree two branched covers of $\Sigma$ that are not branched over any pair of adjacent rays, and they fall into three combinatorial types, as follows.  The fan $\Sigma$ is combinatorially equivalent to the fan whose maximal cones are the cones over the faces of the cube with vertices $\{ (\pm1, \pm1, \pm1) \}$.  (One may think of taking the ray through $(1,1,1)$ and rotating it to the ray through $(1,2,3)$, deforming the fan over the faces of the cube along the way, to get to $\Sigma$.)  The rays of $\Sigma$ correspond to the vertices of the cube, so the rays over which $\bf$ is ramified correspond to a subset of the vertices of the cube, and we may assume that no two of these vertices span an edge.  We represent such a branched cover by drawing a cube with this subset of vertices marked.  Up to automorphisms of the cube, there are three possibilities, shown below.

\begin{center}
\begin{picture}(305,100)(0,-20)

\put(0,0){\line(0,1){50}}
\put(0,0){\line(1,0){50}}
\put(50,0){\line(0,1){50}}
\put(0,50){\line(1,0){50}}
\put(50,0){\line(1,1){15}}
\put(0,50){\line(1,1){15}}
\put(65,15){\line(0,1){50}}
\put(15,65){\line(1,0){50}}
\put(50,50){\line(1,1){15}}

\put(0,0){\circle*{4}}
\put(65,65){\circle*{4}}

\put(30,-20){\makebox(0,0)[b]{(A)}}

\put(120,0){\line(0,1){50}}
\put(120,0){\line(1,0){50}}
\put(170,0){\line(0,1){50}}
\put(120,50){\line(1,0){50}}
\put(170,0){\line(1,1){15}}
\put(120,50){\line(1,1){15}}
\put(185,15){\line(0,1){50}}
\put(135,65){\line(1,0){50}}
\put(170,50){\line(1,1){15}}

\put(120,50){\circle*{4}}
\put(185,65){\circle*{4}}

\put(150,-20){\makebox(0,0)[b]{(B)}}

\put(240,0){\line(0,1){50}}
\put(240,0){\line(1,0){50}}
\put(290,0){\line(0,1){50}}
\put(240,50){\line(1,0){50}}
\put(290,0){\line(1,1){15}}
\put(240,50){\line(1,1){15}}
\put(305,15){\line(0,1){50}}
\put(255,65){\line(1,0){50}}
\put(290,50){\line(1,1){15}}

\put(240,0){\circle*{4}}
\put(290,50){\circle*{4}}
\put(305,15){\circle*{4}}
\put(255,65){\circle*{4}}

\put(270,-20){\makebox(0,0)[b]{(C)}}

\end{picture}
\end{center}

\noindent There are four branched covers of type (A), twelve of type (B), and two of type (C).  To prove the proposition, it remains to show that there are no nontrivial piecewise-linear functions on each of these eighteen branched covers.  These are straightforward linear algebra computations that can be done by hand (or with a computer, if preferred).  We illustrate with the computations for one of these branched covers.

Suppose $\bf$ is ramified exactly over $\rho_1$, $\rho_3$, $\rho_6$, and $\rho_8$, where $\rho_i$ is the ray spanned by $v_i$.  This is one of the two cases of type (C).  We can label the rays of $\Delta$
\[
\{\mbox{rays of $\Delta$} \} = \{ \rho_{11}, \rho_{21}, \rho_{22}, \rho_{31}, \rho_{41}, \rho_{42}, \rho_{51}, \rho_{52}, \rho_{61}, \rho_{71}, \rho_{72}, \rho_{81} \},
\]
such that $f(\rho_{ij}) = \rho_i$, and the maximal cones of $\Delta$ are
\[
\begin{array}{lll}
\sigma_{11} = \< v_{11}, v_{21}, v_{31}, v_{42} \>, & \sigma_{12} = \< v_{11}, v_{22}, v_{31}, v_{41} \>, & \sigma_{21} = \< v_{51}, v_{61}, v_{72}, v_{81} \>, \\
\sigma_{22} = \< v_{52}, v_{61}, v_{71}, v_{81} \>, & \sigma_{31} = \< v_{11}, v_{21}, v_{52}, v_{61} \>, & \sigma_{32} = \< v_{11}, v_{22}, v_{51}, v_{61} \>, \\
\sigma_{41} = \< v_{21}, v_{31}, v_{61}, v_{72} \>, & \sigma_{42} = \< v_{22}, v_{31}, v_{61}, v_{71} \>, & \sigma_{51} = \< v_{31}, v_{41}, v_{72}, v_{81} \>, \\
\sigma_{52} = \< v_{31}, v_{42}, v_{71}, v_{81} \>, & \sigma_{61} = \< v_{11}, v_{41}, v_{52}, v_{81} \>, & \sigma_{62} = \< v_{11}, v_{42}, v_{51}, v_{81} \>,
\end{array}
\]
where $v_{ij}$ is the primitive generator of $\rho_{ij}$.

A piecewise-linear function $\Psi \in \PL(\DDelta) \otimes_\Z \R$ is determined by the values $z_{ij} = \Psi(v_{ij})$, and these $z_{ij}$ must satisfy certain relations to ensure that they come from a linear function on each maximal cone $\sigma_{ij}$.  For instance, since $-2v_1 + 4v_2 -3v_3 + 5v_4 = 0$, if $z_{11}$,  $z_{21}$, $z_{31}$, and $z_{42}$ come from a linear function on $\sigma_{11}$ then
\[
-2z_{11} + 4z_{21} -3z_{31} + 5z_{42} = 0.
\]
Each of the other eleven maximal cones of $\DDelta$ gives a similar linear relation on the $\{ z_{ij} \}$, and it follows that $\Psi \mapsto \{ \Psi(v_{ij}) \}$ takes $\PL(\DDelta) \otimes_\Z \R$ isomorphically onto the set of $\{ z_{ij} \}$ such that
\[
\left(
\begin{array}{rrrrrrrrrrrr}
-2 & 4 & 0 & -3 & 0 & 5 & 0 & 0 & 0 & 0 & 0 & 0 \\
-2 & 0 & 4 & -3 & 5 & 0 & 0 & 0 & 0 & 0 & 0 & 0 \\
0 & 0 & 0 & 0 & 0 & 0 & -1 & 0 & 1 & 0 & -1 & 1 \\
0 & 0 & 0 & 0 & 0 & 0 & 0 & -1 & 1 & -1 & 0 & 1 \\
2 & -4 & 0 & 0 & 0 & 0 & 0 & -3 & 5 & 0 & 0 & 0 \\
2 & 0 & -4 & 0 & 0 & 0 & -3 & 0 & 5 & 0 & 0 & 0 \\
0 & 1 & 0 & -1 & 0 & 0 & 0 & 0 & -1 & 0 & 1 & 0 \\
0 & 0 & 1 & -1 & 0 & 0 & 0 & 0 & -1 & 1 & 0 & 0 \\
0 & 0 & 0 & 1 & -1 & 0 & 0 & 0 & 0 & 0 & -1 & 1 \\
0 & 0 & 0 & 1 & 0 & -1 & 0 & 0 & 0 & -1 & 0 & 1 \\
2 & 0 & 0 & 0 & -5 & 0 & 0 & -3 & 0 & 0 & 0 & 4 \\
2 & 0 & 0 & 0 & 0 & -5 & -3 & 0 & 0 & 0 & 0 & 4 
\end{array}
\right)
\left(
\begin{array}{c}
z_{11} \\
z_{21} \\
z_{22} \\
z_{31} \\
z_{41} \\
z_{42} \\
z_{51} \\
z_{52} \\
z_{61} \\
z_{71} \\
z_{72} \\
z_{81}
\end{array} 
\right)
= 0
\]
This $ 12 \times 12 $ matrix has rank nine, so the solution space is three-dimensional.  Since $\bf^* M_\R$ gives a three-dimensional space of solutions, it follows that every piecewise-linear function on $\DDelta$ is the pullback of a linear function on $N_\R$, and hence is trivial.  Similar computations show that every piecewise-linear function on each of the other seventeen branched covers is the pullback of a linear function on $N_\R$, and it follows that every rank two toric vector bundle on $X$ is trivial.

However, there do exist nontrivial toric vector bundles of rank three on $X$.  To see this, let $E$ be a three-dimensional $k$-vector space with distinct one-dimensional subspaces $L$, and $L'$, and distinct two-dimensional subspaces $V$, $V'$, $V''$, and $V'''$, satisfying the relations
\[
\begin{array}{ccccc}
L \subset V, & L' \subset V', & L' \subset V''', & L \not \subset V'', & L' \not \subset V'',
\end{array}
\]
and
\[
\dim (V \cap V' \cap V'') = \dim (V \cap V' \cap (L + L')) = \dim (V \cap V'' \cap (L + L')) = 1.
\]
Consider the following filtrations.
\[
\begin{array}{ll}
{ E^{\rho_1}(i) = \left\{ \begin{array}{ll} E & \mbox{ for } i  \leq -1,  \\ V & \mbox{ for } i=0, \\ L & \mbox{ for } i  = 1, \\ 0 & \mbox{ for } i > 1, \end{array} \right. } &
{ E^{\rho_2}(i) = \left\{ \begin{array}{ll} E & \mbox{ for } i  \leq 0,  \\ L+ L' & \mbox{ for } 0  <  i   \leq 2, \\ 0 & \mbox{ for } i  > 2, \end{array} \right. } \\ \\ 
{ E^{\rho_3}(i) = \left\{ \begin{array}{ll} E & \mbox{ for } i  \leq -1,  \\ L' & \mbox{ for } 0 < i \leq 2, \\  0 & \mbox{ for } i > 2, \end{array} \right. } &
{ E^{\rho_4}(i) = \left\{ \begin{array}{ll} E & \mbox{ for } i  \leq -2,  \\ V' & \mbox{ for } -2  <  i   \leq 0, \\ 0 & \mbox{ for } i  > 0, \end{array} \right.  } \\ \\
{ E^{\rho_5}(i) = \left\{ \begin{array}{ll} E & \mbox{ for } i  \leq -2,  \\ V'' & \mbox{ for } -2 < i \leq 0, \\  0 & \mbox{ for } i > 0, \end{array} \right. } &
{ E^{\rho_6}(i) = \left\{ \begin{array}{ll} E & \mbox{ for } i  \leq 0,  \\ L' & \mbox{ for } 0  <  i   \leq 2, \\ 0 & \mbox{ for } i  > 2, \end{array} \right. } \\ \\ 
{ E^{\rho_7}(i) = \left\{ \begin{array}{ll} E & \mbox{ for } i  \leq -2,  \\ V''' & \mbox{ for } -2 \leq i \leq 0, \\ L' & \mbox{ for } 0 < i \leq 2, \\ 0 & \mbox{ for } i > 2, \end{array} \right. } &
{ E^{\rho_8}(i) = \left\{ \begin{array}{ll} E & \mbox{ for } i  \leq -2,  \\ L' & \mbox{ for } -2  <  i   \leq 0, \\ 0 & \mbox{ for } i  > 0, \end{array} \right. } \\ \\ 
\end{array}
\]
Then it is straightforward to check that these $\{E^\rho(i)\}$ satisfy the compatibility condition in Klyachko's Classification Theorem, for 
the following multisets of linear functions $\bu(\sigma_i) \subset M$, giving a rank three toric vector bundle on $X$ with nontrivial Chern class.  The multisets are
\[
\begin{array}{l}
\bu(\sigma_1) = \bigl\{ (1,-1,0), (0,-1,1), (0,0,0) \bigr\}, \vspace{2 pt} \\
\bu(\sigma_2) = \bigl\{ (0,-1,1), (0,-1,-1), (1,0,1) \bigr\}, \vspace{2 pt}\\ 
\bu(\sigma_3) = \bigl\{ (1,-1,0), (0,-1,1), (0,0,0) \bigr\}, \vspace{2 pt} \\
\bu(\sigma_4) = \bigl\{ (1,0,1), (0,-2,0), (0,0,0) \bigr\}, \vspace{2 pt}\\
\bu(\sigma_5) = \bigl\{ (1,-1,0), (-1,-1,0), (1,0,1) \bigr\}, \vspace{2 pt}\\ 
\bu(\sigma_6) = \bigl\{ (1,-1,0), (0,-1,1), (0,0,0) \bigr\}.
\end{array}
\]
\end{example}

\begin{remk}
The multisets $\bu(\sigma_i)$ associated to the nontrivial rank three toric vector bundle in Example \ref{Fulton's threefold} were found by a computer program searching through the groups of piecewise-linear functions on all of the maximal degree three branched covers of the given fan.  With such a piecewise-linear function in hand, in such a relatively simple situation, it is not difficult to construct a family of filtrations $\{ E^\rho(i) \}$ giving a nontrivial toric vector bundle (in spite of the negative results in \cite[Section~4]{moduli}).   
\end{remk}

\begin{example} \label{new example}
Let $N = \Z^3$, and let $\Sigma'$ be the fan in $N_\R$ combinatorially equivalent to the fan $\Sigma$ in Example \ref{Fulton's threefold} that is obtained by replacing the ray through $(1, -1, 1)$ with the ray through $(1, -1, 2)$.  In other words, $\Sigma'$ is the fan whose rays are spanned by
\[
\begin{array}{llll}
v'_1 = (1,2,3), & v'_2 = (1,-1,2), & v'_3 = (-1,-1,1), & v'_4 = (-1,1,1), \\
v'_5 = (1,1,-1), & v'_6 = (1,-1,-1), & v'_7 = (-1,-1,-1), & v'_8 = (-1,1,-1),
\end{array}
\]
and whose maximal cones are
\[
\begin{array}{lll}
\sigma'_1 = \< v'_1,v'_2, v'_3, v'_4\>, & \sigma'_2 = \<v'_5, v'_6, v'_7, v'_8\>, & \sigma'_3 = \< v'_1, v'_2, v'_5, v'_6 \>, \\
\sigma'_4 = \< v'_2, v'_3, v'_6, v'_7\>, & \sigma'_5 = \< v'_3, v'_4, v'_7, v'_8 \>, & \sigma_6 = \<v'_1, v'_4, v'_5, v'_8 \>.
\end{array}
\]
The toric variety $X' = X(\Sigma')$ is complete, and computations similar to those in Example \ref{Fulton's threefold} show that $X'$ has no nontrivial toric vector bundles of rank less than or equal to two.

We claim that $X'$ has no nontrivial rank three toric vector bundles.  Since $X'$ has no nontrivial line bundles, there are no nonconstant morphisms from $X'$ to projective varieties, so Theorem \ref{trivial rank 3} implies that any rank three toric vector bundle on $X'$ with trivial Chern class is trivial.  Therefore, to prove the claim it suffices to show that there are no nontrivial piecewise-linear functions on any maximal degree three branched cover of $\Sigma'$.

Let $B'$ be the finite set of interesection points of the rays of $\Sigma'$ with the unit sphere $S$ in $\R^3 = N_\R$.  By Theorem \ref{bijections}, the number of isomorphism classes of maximal branched covers of $\Sigma'$ is equal to the number of conjugacy classes of representations
\[
\pi_1(S \smallsetminus B') \rightarrow S_3.
\]
Since $\pi_1(S \smallsetminus B')$ is free on seven generators, there are $6^7 = 279,936$ such representations.  By considering conjugacy classes, and by excluding certain patterns of monodromy around points corresponding to adjacent rays, as in the analysis of piecewise-linear functions on degree two branched covers in Example \ref{Fulton's threefold}, one could reduce the number of cases to consider by a couple of orders of magnitude, but there would still be too many to comfortably check by hand.  To verify the claim, we programmed a computer to step through each of the $279,936$ possible monodromy representations, computing generators for the group of piecewise-linear functions on the corresponding branched cover in each case, and checking whether each of the generators is trivial.  The program that was used consists of roughly one hundred lines of code in \emph{Mathematica}, runs on a moderately fast desktop computer (by 2005 standards) in approximately one hour, and confirms that there are no nontrivial piecewise-linear functions on any maximal degree three branched covers of $\Sigma'$.  The source code is available from the author on request.  
\end{example}

\bibliography{math}
\bibliographystyle{amsalpha}

\end{document}